\documentclass[12pt,draftcls,onecolumn]{IEEEtran}
\ifCLASSINFOpdf
  \usepackage[pdftex]{graphicx}
	
\else

\fi
\usepackage{mathrsfs}
\usepackage{amsthm}
\usepackage{amsmath}
\usepackage{amsfonts}
\usepackage{epstopdf}
\usepackage{tabularx}
\usepackage{color}
\usepackage{mathtools}
\usepackage[utf8]{inputenc}
\usepackage[T1]{fontenc}
\usepackage{array}
\usepackage{multirow}
\usepackage[]{algorithm}
\usepackage[]{algpseudocode}
\usepackage[]{pseudocode}
\usepackage[caption=false,font=footnotesize]{subfig}

\DeclareMathOperator*{\argmin}{arg\,min}
\DeclareMathOperator*{\polylog}{polylog}
\usepackage{color}

\newcommand{\cA}{\mathcal{A}}
\DeclarePairedDelimiter{\dotp}{\langle}{\rangle}
\begin{document}
\normalsize
%
\title{Beating level-set methods for 3D seismic data interpolation: \\a primal-dual alternating approach}

\author{Rajiv Kumar$^1$, Oscar L\'{o}pez$^2$, Damek Davis$^{^3}$, Aleksandr Y. Aravkin$^{^4}$ and Felix J. Herrmann$^1$\\
$^1$Department of Earth, Ocean and Atmospheric Sciences, The University of British Columbia \\
$^2$Department of Mathematics, The University of British Columbia\\
$^3$ Department of Operations Research and Information Engineering, Cornell University\\
$^4$ Department of Applied Mathematics, University of Washington}

\maketitle

\normalsize
\begin{abstract}
Acquisition cost is a crucial bottleneck for seismic workflows, and 
low-rank formulations for data interpolation allow practitioners to `fill in'  data volumes from critically subsampled data  
acquired in the field. Tremendous size of seismic data volumes required for seismic processing
remains a major challenge for these techniques. 

We propose a new approach to solve residual constrained formulations 
for interpolation. We represent the data volume using matrix factors, and 
build a block-coordinate algorithm with constrained convex subproblems 
that are solved with a primal-dual splitting scheme. 
The new approach is competitive with state of the art level-set algorithms that interchange the role 
of objectives with constraints. We use the new algorithm 
to successfully interpolate a large scale  5D seismic data volume,
generated from the geologically complex  synthetic 3D Compass velocity model, 
where 80\% of the data has been removed. 

\end{abstract}


\section{Introduction}
\label{sec:intro}
Seismic data interpolation is crucial for accurate inversion and imaging procedures such as full-waveform inversion \cite{virieux2009overview}, reverse-time migration \cite{baysal1983reverse, Plessix} and multiple removal methods (SRME, EPSI) \cite{verschuur1992adaptive, lin2013GEOPrepsi}.
Dense acquisition is prohibitively expensive in these applications, motivating reduction in seismic measurements. 
On the other hand, using subsampled sources and receivers without interpolation gives unwanted imaging artifacts. 
A range of trace interpolation methodologies have been proposed, exploiting low dimensional structure of seismic data, including sparsity (\cite{curvelet, yilmaz}) and low-rank (\cite{slim,kreimer2013tensor,oropeza2011simultaneous,trickett2010rank,osher}), and with theoretical guarantees 
for low-rank matrix recovery available in a range of contexts~(\cite{fazel,recht}). 
The main goal is to simultaneously sample and compress a signal using optimization to replace dense acquisition, 
thus enabling a range of applications in seismic data processing at a fraction of the cost.

Low-rank matrix completion techniques (\cite{candesrecht,candes}) have been successfully applied
 to seismic trace interpolation. 
Here we focus on {\it residual-constrained} formulations, which minimize a regularizer subject to 
a constraint on the data misfit. 
This formulation is particularly appealing when practitioners have an estimate of the {\it noise floor}~\cite{spgl1,slim}. 

In \cite{slim}, the authors combine explicit low-rank factorization~(\cite{recht,srebro}) with {\it level-set optimization techniques}~(\cite{spgl1,AravkinBurkeFriedlander:2013,aravkin2016level}), and apply the resulting approach 
to seismic data interpolation in the midpoint-offset domain. A recent extension~\cite{OffGrid} includes seismic data sampled on unstructured grids, with recovery error bounds provided for the methodology.

Here, we develop a competitive optimization scheme based on alternating minimization for factorized 
formulations. The alternating approach in our context yields
convex residual-constrained subproblems,  which we solve using primal-dual splitting techniques of~\cite{chambolle2011first}. 
The resulting scheme is simple, and can be applied to extremely large-scale problems.

The paper proceeds as follows. In Section~\ref{sec:LR}, we formulate seismic data interpolation as a 
low-rank optimization problem, highlighting the transform domain that makes low-rank  techniques applicable. 
In Section~\ref{sec:meth}, we present the alternating optimization scheme for the residual-constrained formulation, 
together with an algorithm to solve the convex subproblems in each factor. The new approach is 
illustrated using 3D seismic data volumes in Section~\ref{sec:exp}, and we show it compares favorably
to the the level-set approach of~\cite{slim}. We end with a discussion and future directions in Section~\ref{sec:discussion}.

\section{Notation}
\label{sec:notation}
We use lower case to represent vectors ($b,f$), upper case to denote matrices and tensors ($X, Y$) , and calligraphic upper case  for linear operators ($\mathcal{A}$). 
For 3D seismic data acquisition, seismic data volumes have two source dimensions ($s_x, s_y$), two receiver dimensions ($r_x, r_y$), and time $t$.

\section{Low-rank matrix completion}
\label{sec:LR}
The goal of low-rank matrix completion is to accurately estimate the unobserved entries of a data matrix, $X\in\mathbb{C}^{n\times m}$, from the observed entries and prior knowledge that the matrix exhibits low-rank structure, i.e. has few nonzero singular values, or can be accurately approximated by such a low-rank matrix. Letting $\Omega\subset\{1,2,...,n\}\times\{1,2,...,m\}$ be the set of observed entries, we  define a sampling operator $P_{\Omega}:\mathbb{C}^{n\times m}\mapsto\mathbb{C}^{n\times m}$ by its element-wise action:
\[
    P_{\Omega}(X)_{ij}=\left\{
                \begin{array}{ll}
                  X_{ij}, \ \ \mbox{if} \ \ (i,j)\in\Omega,\\
                  0, \ \ \mbox{otherwise}.
                \end{array}
              \right.
 \]
We can write our observations as $b = P_{\Omega}(X) + \epsilon$, where $\epsilon = P_{\Omega}(\epsilon)\in\mathbb{C}^{n\times m}$ models data corruption on the subset of observed entries, with a given noise floor $\|\epsilon\|_F\leq\eta$. 
Rather than minimizing the rank of $X$, which is a combinatorial problem, 
we can consider the nuclear norm $\|X\|_*$; in many contexts this allows exact recovery guarantees on $X$~\cite{fazel,recht}.
This convex relaxation gives us the interpolated data volume $X^{\sharp}$ minimizer of 
\begin{equation}
X^{\sharp} := \argmin_{X\in \mathbb{C}^{n\times m}} \|X\|_* \ \ \mbox{s.t.} \ \ \|P_{\Omega}(X)-b\|_F \leq \eta,
\label{eq:NN}
\end{equation}
where $\|\cdot\|_F$ denotes the Frobenius norm. This procedure finds the $X$ with lowest 
nuclear norm that still fits observations up to the noise level $\eta$.
In real-world applications, $\eta$ may require estimation by cross-validation or other techniques. 

Matrix completion via nuclear norm minimization has been extensively studied (\cite{recht,candesrecht,candes,fazel}). Typical results show that if $\Omega$ is generated uniformly randomly, then any $n\times m$ ``incoherent'' rank $r$ matrix can be reconstructed with as few as $|\Omega| \sim\mathcal{O}(\max(n,m)r\polylog{\max(n,m)})$ observed entries. When $r\ll \min(n,m)$, this methodology provides a great reduction in the number of measurements needed compared to typical sampling procedures used in many applications. These insights have spurred work with the goal of solving (\ref{eq:NN}) efficiently.

Many implementations to solve (\ref{eq:NN}), such as singular-value projection and singular value thresholding (\cite{jain,SVT}) 
require singular value decomposition (SVD) or partial SVD computations of the data matrix. 
Furthermore, treating the entire $X$ as a decision variable requires storage and manipulation of an $n\times m$ object 
at every iteration. Storage and manipulation of large volumes becomes prohibitive for huge-scale problems, making (\ref{eq:NN}) and equivalent formulations impractical for realistic seismic interpolations where the spatial interpolated grid for 3D seismic data acquisition is of the order $10^9$. To rectify these issues (\cite{slim, recht, srebro}) propose a matrix factorization approach. Stipulating a maximal rank $r$ for $X$, 
the authors write the representation $X = LR^{H}$, with $L\in \mathbb{C}^{n\times r}$, 
$R\in \mathbb{C}^{m\times r}$ and $^{*}$ denoting the Hermitian transpose. 
As shown in (\cite{recht,srebro}),
\[
\|X\|_* = \min_{\{L,R: X = LR^{H}\}}\frac{1}{2}(\|L\|_F^2 + \|R\|_F^2)  = \min_{\{L,R: X = LR^{H}\}}\frac{1}{2}\left\| \begin{bmatrix} L \\ R \end{bmatrix}\right\|_F^2,
\]
motivating the formulation of~\cite{slim}: 
\begin{equation}
\min_{L\in \mathbb{C}^{n\times r}, R\in \mathbb{C}^{m\times r}} \frac{1}{2}\left\| \begin{bmatrix} L \\ R \end{bmatrix}\right\|_F^2 \ \ \mbox{s.t.} \ \ \|P_{\Omega}(LR^{H})-b\|_F \leq \eta.
\label{eq:LR}
\end{equation}
If $r\ll\min(n,m)$, we have reduced the memory requirements from $mn$ to $rn + rm$, 
while projection onto the Frobenius norm ball is linear in the size of the variable, 
avoiding expensive SVD computations. 

When $X$ is a full-rank matrix that can be accurately approximated by a rank $r$ matrix, 
implementing (\ref{eq:LR}) with this choice of $r$ provides reconstruction error bounds proportional to the error of the best $r$-rank approximation of $X$ (see for example \cite{candes}). 
This holds true for many applications, where e.g. the data matrix is not low-rank but exhibits quickly decaying singular values, so that the error of the best $r$-rank approximation will be small for appropriately chosen $r$. 
When the underlying rank is not known, minimizing a regularization functional subject to a data constraint~\eqref{eq:LR} has an important practical consequence: as the nominal rank $r$ (number of columns in $L$ and $R$) increases, 
we do not overfit the data~\cite{slim}. 

Simply considering standard seismic data volumes does not immediately yield low-rank structure. 
In the next section, we show how seismic data can be transformed in order to find and  
exploit low-rank representation for interpolation.\\

\subsection{Low-rank structure of seismic data}

\begin{figure}[htp]
\centering
\subfloat[]{
\includegraphics[width=.24\textwidth]{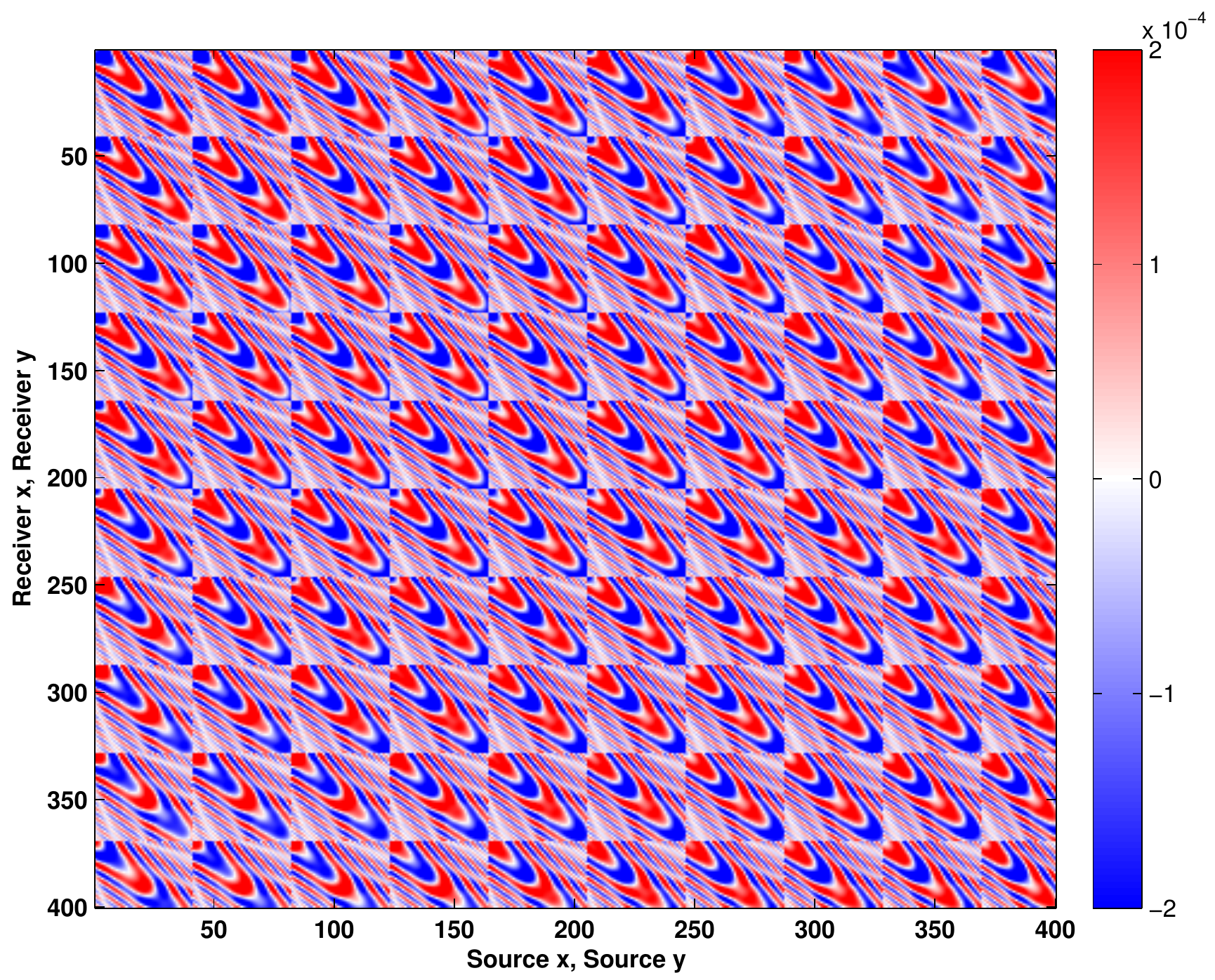}}~
\subfloat[]{
\includegraphics[width=.24\textwidth]{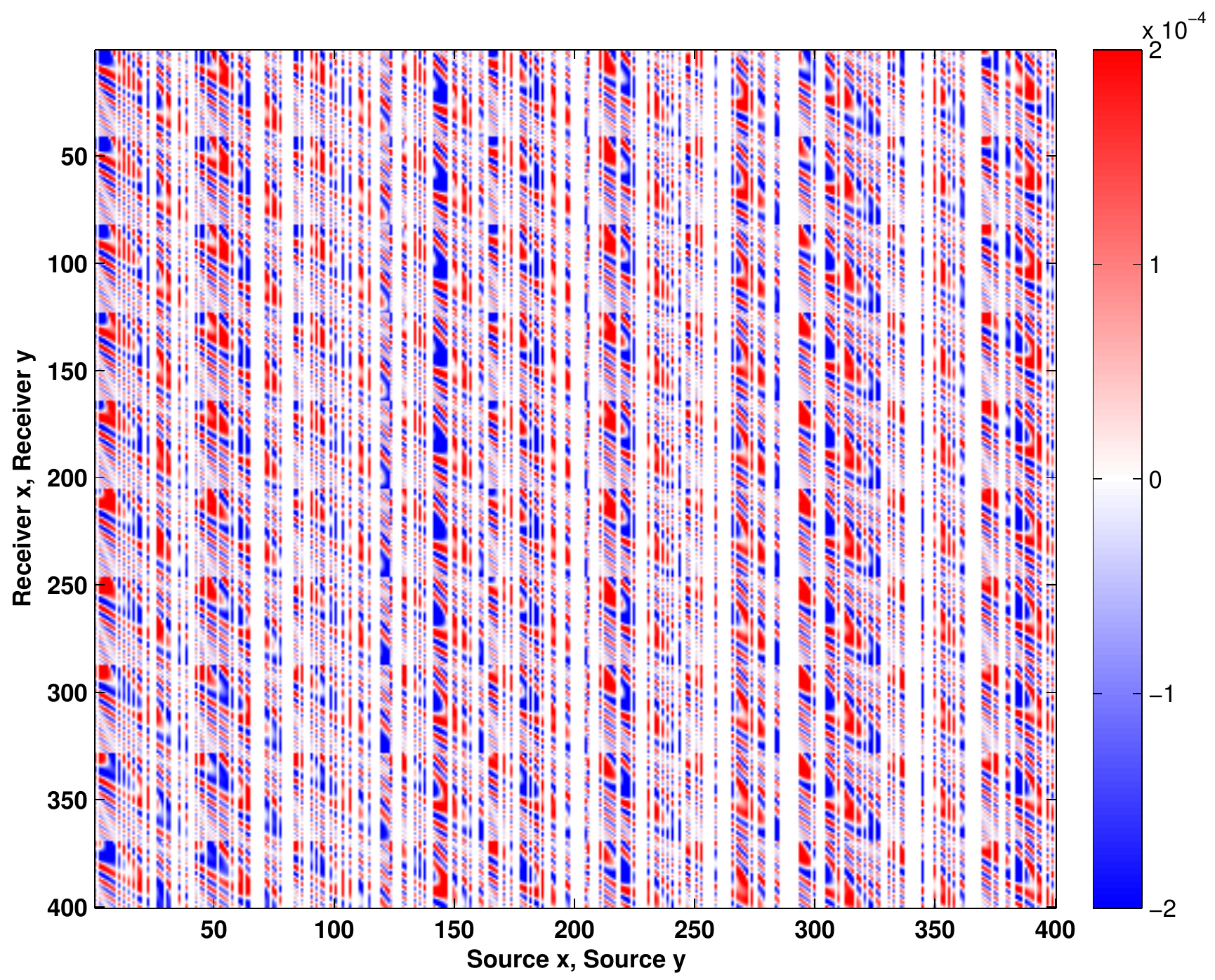}} \\
\subfloat[]{
\includegraphics[width=.24\textwidth]{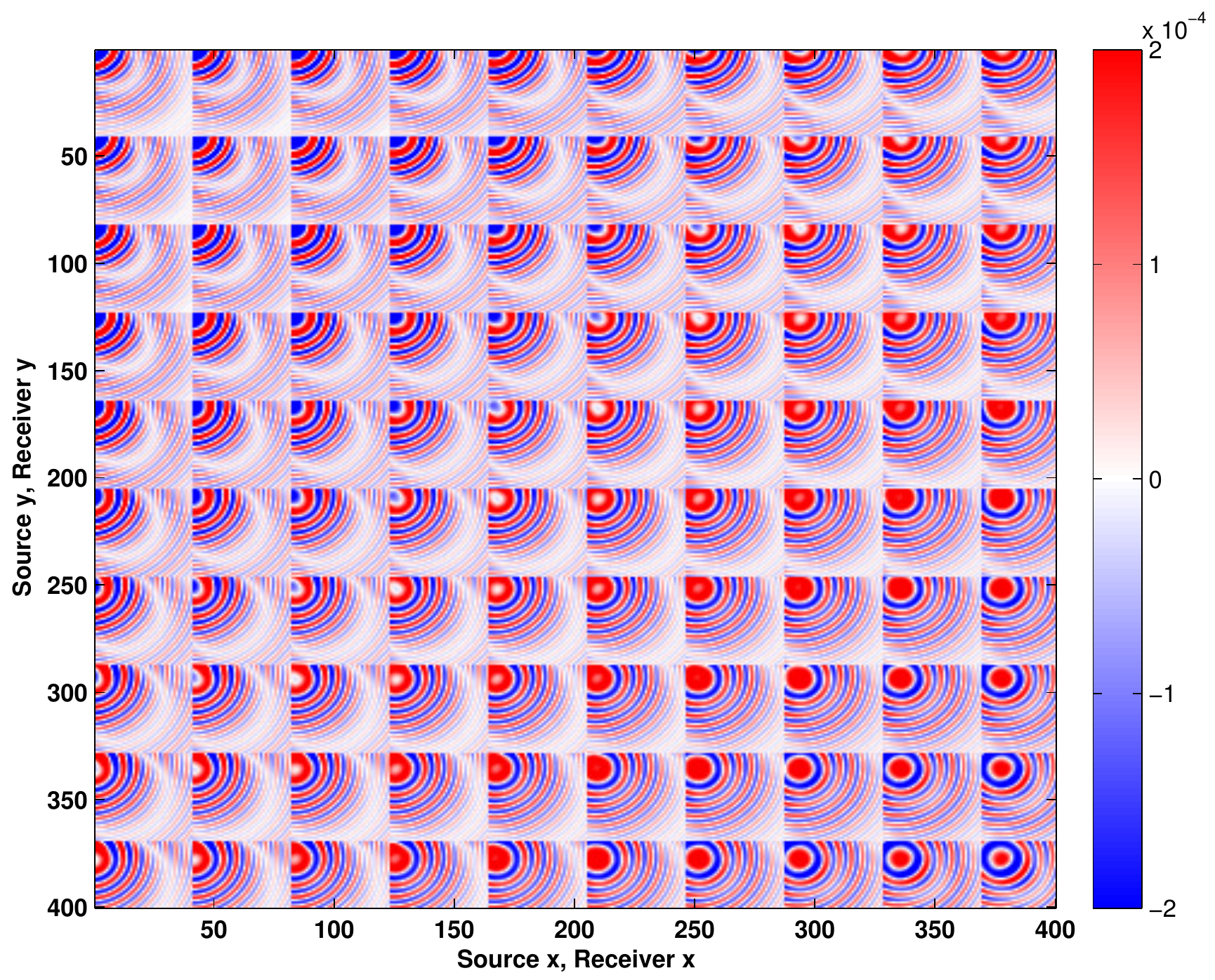}}~
\subfloat[]{
\includegraphics[width=.24\textwidth]{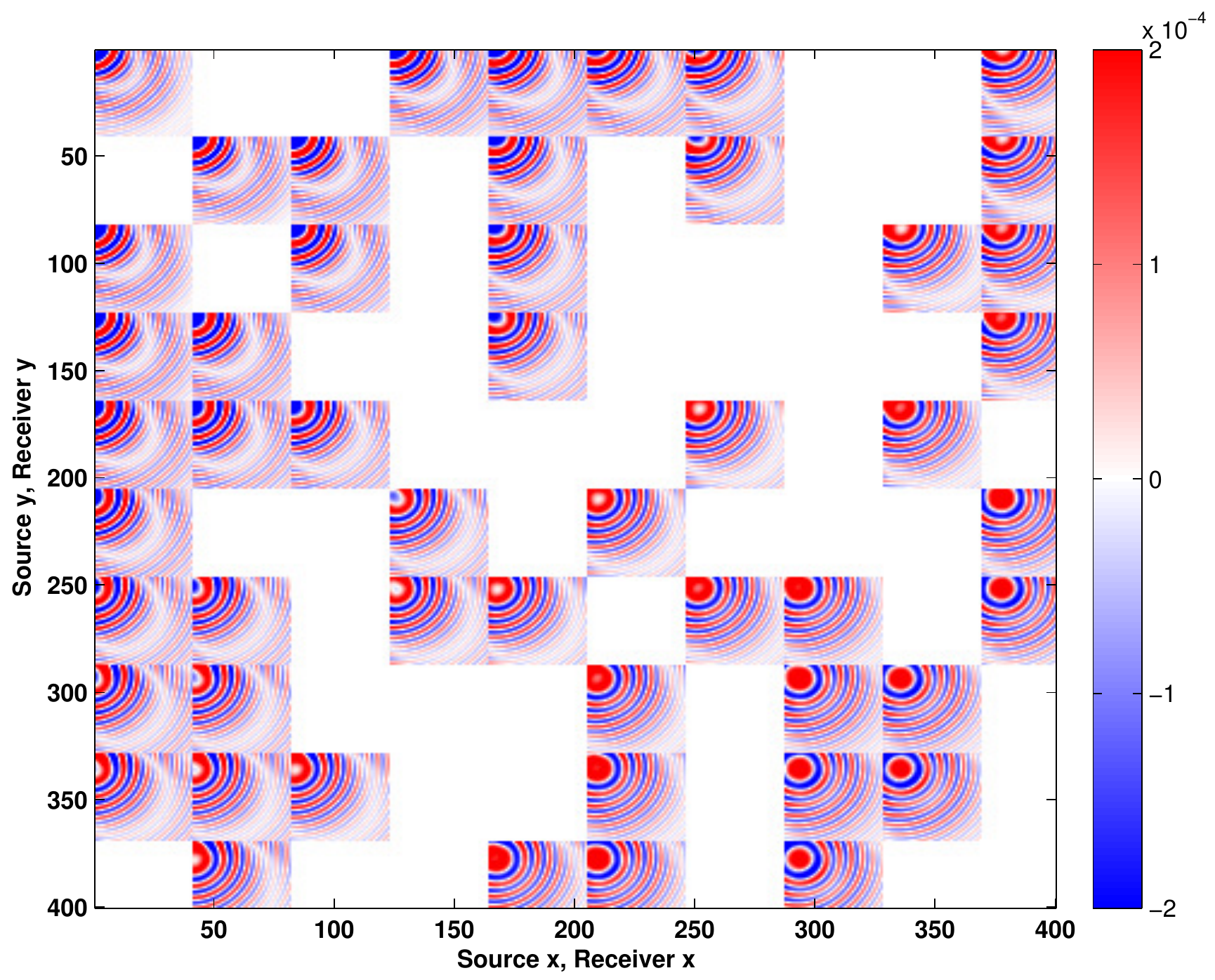}}
\caption{Matricization of 4D monochromatic slice. Top: ($s_x, s_y$) matricization. Bottom: ($r_x, s_x$) matricization.
                  Left: Fully sampled data; Right: Subsampled data.}
\label{fig:3Dmat}
\end{figure}

\begin{figure}[htp]
\centering
\subfloat[]{
\includegraphics[width=.24\textwidth]{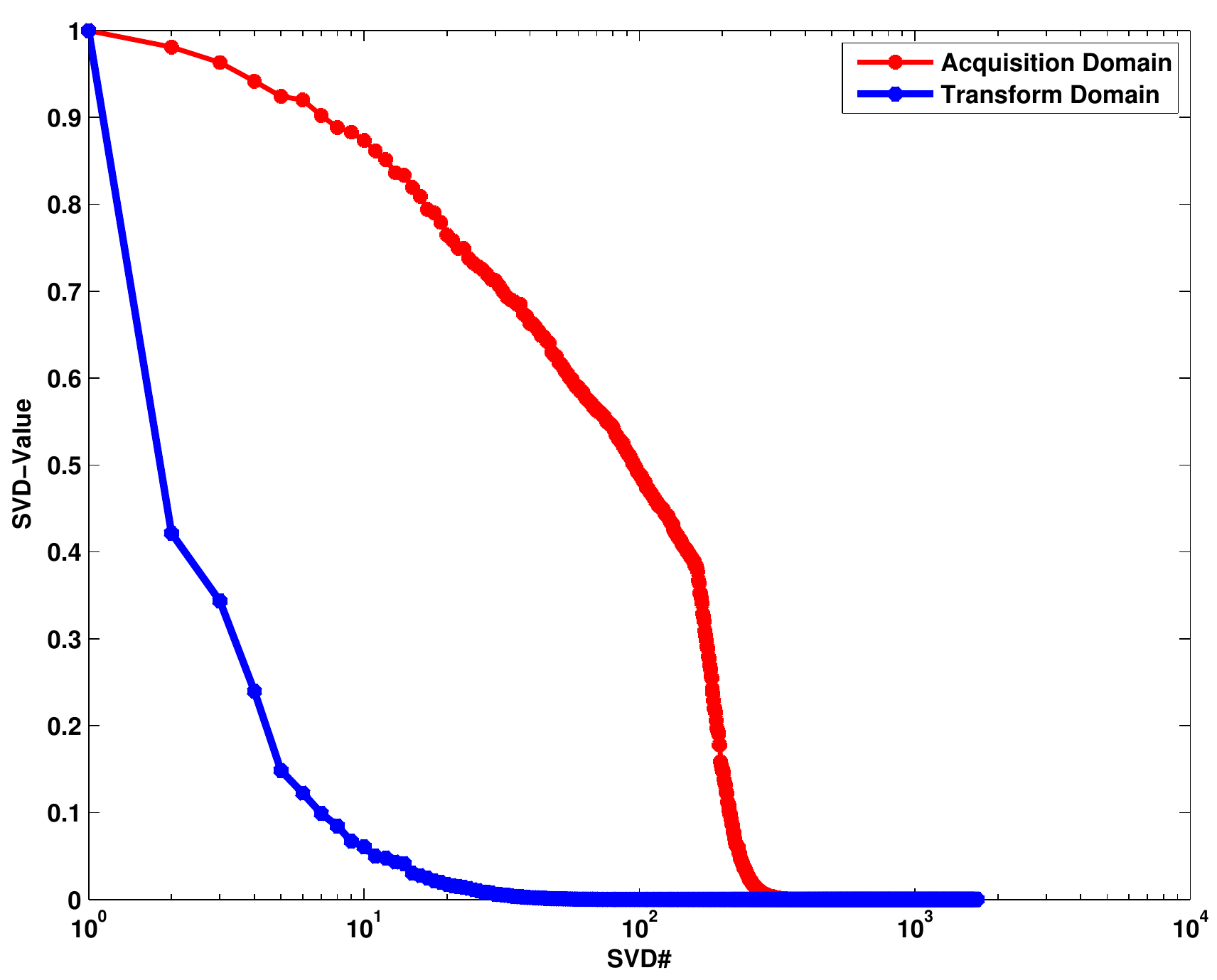}}~
\subfloat[]{
\includegraphics[width=.24\textwidth]{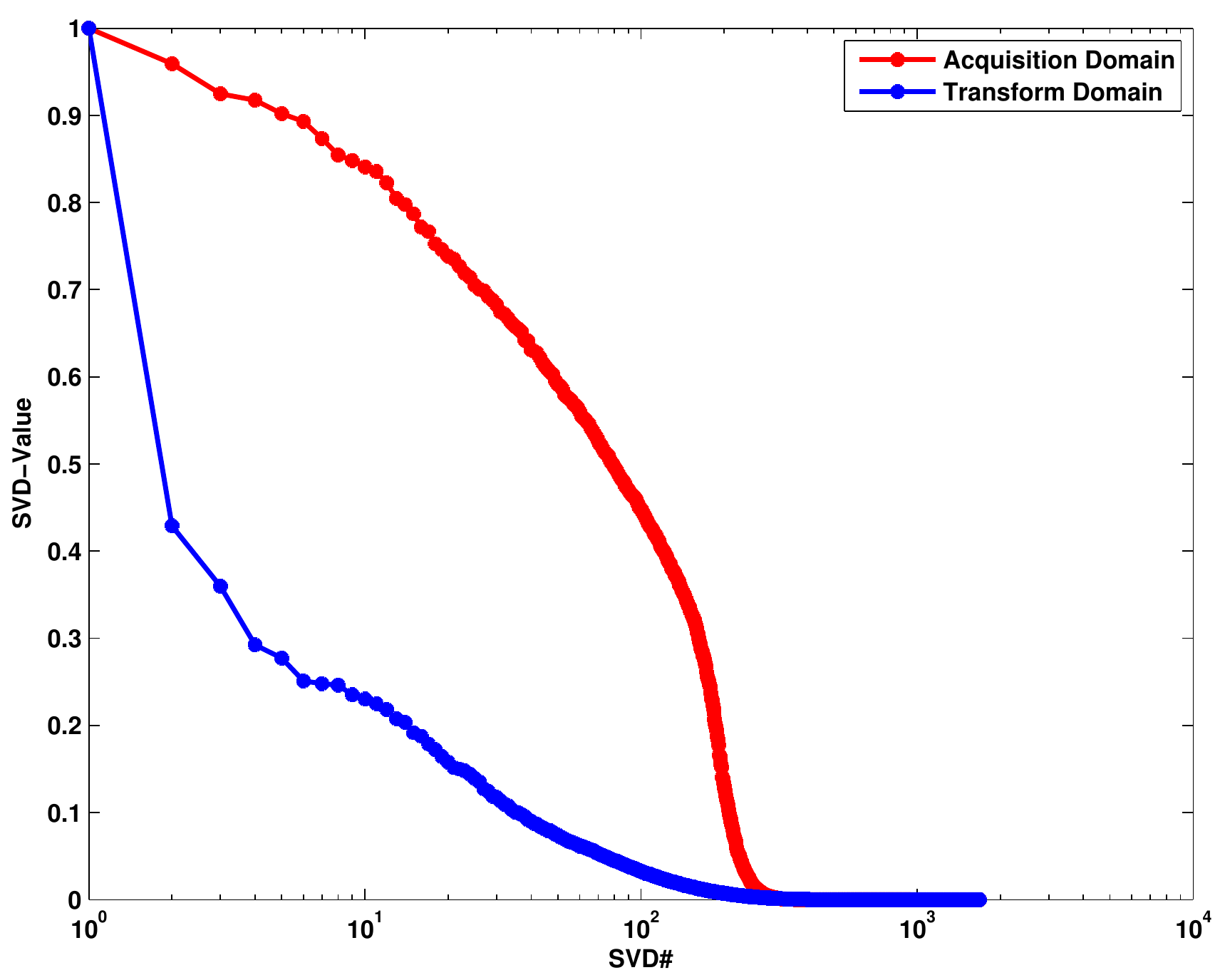}}
\caption{Singular value decay of ($s_x, s_y$) matricization (red curve) and ($r_x, s_x$) matricization (blue curve). (a) Fully sampled data and (b) subsampled data.}
\label{fig:3Dmatsvd}
\end{figure}

For 3D seismic data acquisition, each monochromatic tensor is a 4-dimensional volume with source dimensions $s_x, s_y$, and receiver dimensions $r_x$, $r_y$, respectively. In order to exploit low-rank structure for 3D seismic data acquisition, \cite{dasilva2014htuck, trans} proposed two choices of {\it matricization}, or operation of unfolding a tensor into a matrix along specific dimensions. We either place the $r_x, r_y$ dimensions in the rows and $s_x, s_y$ dimensions in the columns (Figure \ref{fig:3Dmat}a), or $r_y, s_y$ dimensions in the rows and $r_x, s_x$ dimensions in the columns (Figure \ref{fig:3Dmat}c). In $r_x, r_y$ matricization, operator $P_{\Omega}$ samples observed rows and/or columns, i.e. observed data corresponding to source and/or receiver combinations. Placing both receiver coordinates along the rows, i.e., matricization $r_x, r_y$, yields a matrix that has high rank (Figure \ref{fig:3Dmatsvd}a, red curve) and action of the subsampling operator (in particular removal of rows and/ or columns) can only {\it decrease} the rank (Figure \ref{fig:3Dmatsvd}b, red curve). Therefore, interpolation of the data by rank penalization is doomed to failure when using $r_x, r_y$ matricization, since the fully sampled object has higher rank (and larger singular values) than the subsampled dataset. 

We therefore look for a transform domain where the fully sampled data matrix exhibits low-rank structure, while subsampling by $P_\Omega$ increases the rank and/or negatively alters the fast decay of singular values. In this context, we can expect low-rank penalization to help. We use the rank-revealing transforms of \cite{slim, curt4d, trans}. To motivate the use of these transforms, we first explore the singular value decay of full and subsampled volumes in the context of 3D seismic data acquisition. Note that, monochromatic seismic data matrices are full rank even in the transformed domain, but they have fast decay of the singular values in the transform domain, which means that these matrices can be well approximated with low-rank matrices. Therefore, we assume that seismic data matrix exhibit low-rank structure when it has fast-decay of singular values and high-rank structure when it has slow-decay of singular values. In case of 3D seismic data acquisition, matricization $r_x, s_x$ yields fast decay of the singular values (Figure \ref{fig:3Dmatsvd}a, blue curve) for the original fully sampled data volume, while subsampling causes the singular values to decay at a slower rate (Figure \ref{fig:3Dmatsvd}b, blue curve). Therefore, we select the latter matricization for 3D seismic data acquisition.

Working under an appropriate rank-revealing transformation domain, the subsampling operator $P_{\Omega}$ will also involve a transformation operator represented by $\mathcal{S}:\mathbb{C}^{n\times m}\mapsto\mathbb{C}^{p\times q}$. This transformation operator projects the seismic data from the source-receiver domain to the rank-revealing transform domains and its adjoint reverse the operation. We call this sampling-transformation operation
 $\mathcal{A}(\cdot) = P_{\Omega}\mathcal{S}^{H}(\cdot) : \mathbb{C}^{p\times q}\mapsto\mathbb{C}^{n\times m}$, where $^{H}$ represents the adjoint,
 and minimize~\eqref{eq:LR} with $P_{\Omega}$ replaced by $\mathcal{A}$.

\section{Methodology}
\label{sec:meth}
While several formulations are available to solve the equality constrained version of~(\ref{eq:LR}) (i.e. where $\eta = 0$), see~\cite{recht}, and penalized formulations, e.g.~\cite{srebro}, few focus on the case where $\eta$ 
is provided by the user. The level-set factorized approach, which we call LR-BPDN following \cite{slim}, 
solves (\ref{eq:LR}) for a prescribed $\eta$ by defining the value function 
\[
v(\tau) = \min_{\|L\|_F^2 + \|R\|_F^2 \leq 2\tau}  \||\mathcal{A}(LR^{H})-b\|_F,
\]
and applying Newton's method to find $v(\tau) = \eta$. 
In practice 
this means inexactly solving a sequence of optimization problems to evaluate $v(\tau)$, 
and using duality theory to compute $v'(\tau)$ in order to compute the next $\tau$ iterate~\cite{spgl1,slim}. 
The overall approach is reliable but can take many iterations, especially for high-fidelity data fitting.


Here, we propose a block-coordinate descent scheme, outlined in Algorithm \ref{alg:alg1}.

\begin{algorithm}
  \caption{Nuclear Norm Minimization via Alternating Optimization
    \label{alg:alg1}}
  \begin{algorithmic}[1]
		  \Require{$\mathcal{A}$, $b$, $r$, $\eta$, $K_0$}\\
			\textbf{initialize:} any $L^0\in \mathbb{C}^{p\times r}$, $\eta_0 = \|b\|_F$, $\alpha = 0.1$.
			\For{$k = 0 \textrm{ to } K_0-1$}\\
			       $\eta_{k+1} = \max(\alpha^k \eta_{k}, \eta)$.\\
				$\displaystyle R^{k+1} = \argmin_{\{R\in \mathbb{C}^{q\times r}\}} \frac{1}{2}\|R\|_F^2 \ \mbox{s.t.} \  \|\mathcal{A}(L^kR^{H})-b\|_F \leq \eta_{k+1}$\\
				$\displaystyle L^{k+1} = \argmin_{\{L\in \mathbb{C}^{p\times r}\}} \frac{1}{2}\|L\|_F^2 \ \mbox{s.t.} \ \|\mathcal{A}(L{R^{H}}^{k+1})-b\|_F \leq \eta_{k+1}$
			\EndFor\\
			\RETURN{$X^{\sharp} := {\mathcal{S}^{H}}(L^{K_0}\left(R^{H}\right)^{K_0})$}
  \end{algorithmic}
\end{algorithm}

Steps 3 and 4 focus on a single matrix factor at a time. Alternating approaches are 
common for matrix factorization (\cite{altmin,alt,Hastie}). 
The main competing algorithms for constrained matrix completion formulations~\eqref{eq:LR} use a level-set 
approach along the lines of~\cite{slim} (see e.g \cite{altmin}).

We implement steps $3$ and $4$ using a matrix-free approach, described in detail below.

\subsection{Solving Steps 3 and 4 with Primal-Dual Splitting}

All of the work of Algorithm~\ref{alg:alg1} occurs in steps 3 and 4. Each step is an inequality constrained, convex optimization problem. In typical primal-only iterative optimization algorithms, such as the projected-gradient method, we must, at every iteration, project onto the constraint set. In large-scale problems, such a projection itself requires an iterative algorithm. Instead, we solve steps 3 and 4 using the \emph{primal-dual} approach developed in~\cite{chambolle2011first}, each step of which consists solely of matrix-vector products and scalar multiplication, 
and which is completely free of projections, making it suitable for extremely large-scale problems. 
Moreover, it is very simple to implement: the full details of the approach are contained in Algorithm~\ref{alg:alg2} for the $L$-update (step 5).
The $R$-update (step 4) is obtained immediately by changing the roles of $L$ and $R$. 
The matrix norm $\|\cdot\|_{\mathrm{op}}$ used by Algorithm~\ref{alg:alg2} for stepsize selection 
is the largest singular value of the input matrix.  
\begin{algorithm}
  \caption{Primal-Dual Splitting Algorithm for Step 5 in Algorithm \ref{alg:alg1}
    \label{alg:alg2}}
  \begin{algorithmic}[1]
		  \Require{$\mathcal{A}$, $b$, $r$, $\eta$, $K_1$, $R\in \mathbb{C}^{q\times r}$}\\
			\textbf{initialize:} $k=0$, any $y^0\in\mathbb{C}^{n\times m}, L^0\in\mathbb{C}^{p\times r}$ \\
			\quad $\cA : \mathbb{C}^{p\times r} \rightarrow \mathbb{C}^{n\times m}$ defined by $L \mapsto \mathcal{A}(L R^\ast)$ \\
			\qquad $\gamma = \frac{0.99}{\|R\|_{\mathrm{op}}}$ 
			\For{$k = 0 \textrm{ to } K_1-1$}\\
			\quad $L^{k+1} = \frac{1}{1+\gamma}(L^k - \gamma \cA^\ast y^k)$ \\
			\quad $y^+ = y^k + \gamma \cA(2L^{k+1} - L^k) - \gamma b$\\
			\quad $y^{k+1} = \max\left\{1-\frac{\eta\gamma}{\|y^+\|_{F}}, 0\right\} y^+.$
			\EndFor\\
			\RETURN{$L = L^{K_1}$}
  \end{algorithmic}
\end{algorithm} 

To understand Algorithm~\ref{alg:alg2}, 
we formulate a convex-concave saddle-point representation of step~5 of Algorithm~\ref{alg:alg1}:
\begin{equation}
\label{eq:saddle}
\begin{aligned}
\min_L\max_{y} \; \left\{\frac{1}{2}\|L\|^2 + \dotp{\cA L - b, y} - \eta\|y\|_2\right\}.
\end{aligned}
\end{equation}
Computing the maximum over $y$ in~\eqref{eq:saddle} immediately recovers 
the primal problem in $L$ for step~5 of Algorithm~\ref{alg:alg1}, 
since this maximum is equal to $\frac{1}{2}\|L\|^2$ when $\|\cA L - b\| \leq \eta$, 
and is infinity otherwise. 
To align Algorithm~\ref{alg:alg2} with the notation of~\cite{chambolle2011first},
set $G(L) = \frac{1}{2}\|L\|_F^2$ and 
$$
F(z) = \begin{cases}
0 & \text{if $\|z\|_2 \leq \eta$;} \\
\infty & \text{otherwise,}
\end{cases}
$$
i.e., $F$ is the convex indicator of the scaled $\ell_2$ norm ball. 
The convex conjugate $F^\ast$ can be immediately computed: $F^\ast(y) = \eta\|y\|_2$ (see e.g.~\cite{RTR}). 
Then steps~5 and~7 of Algorithm~\ref{alg:alg2} correspond to the \emph{proximal updates}:
\begin{align*}
L^{k+1}  &= \argmin_{L} \left\{G(L) + \frac{1}{2\gamma} \|L - (L^k - \gamma \cA^\ast y^k)\|^2\right\}; \\
y^{k+1}  &= \argmin_{y} \left\{ F^\ast(y) + \frac{1}{2\gamma}\|y - y^+\|^2\right\},
\end{align*}
where $y^+$ is defined as in step~6 of Algorithm~\ref{alg:alg2}. The reader can check that any 
fixed point $(\overline L, \overline y)$ of this iteration solves the saddle point formulation~\eqref{eq:saddle},
and hence the problem in step~5 of Algorithm~\ref{alg:alg1}. 
See~\cite{chambolle2011first} for a detailed convergence analysis and rate results. 

\subsection{Relaxation of the $\eta$-constraint}
\label{sec:relax}
We found that running Algorithm~\ref{alg:alg1} with $\eta_k = \eta$, the final target value, is 
too aggressive, especially in the context of the alternating scheme in $L$ and $R$. In particular, 
in the early steps, when the values of $L$ and $R$ are far from the optimum, it is to our 
advantage to solve approximate subproblems. As the algorithm progresses, 
we tighten $\eta_k$, and eventually solve the problem of interest with $\eta_k = \eta$. 

\section{Numerical Experiments}
\label{sec:exp}
We perform seismic data interpolation using the proposed primal-dual alternating scheme on a 5D synthetic data set from a realistic complex geological model. Our goal is to interpolate 3D geological models with very fine-scale features and complex sedimentary environments (where interpolation problems are very challenging for densely sampled data), comparing level-set and primal-dual methods for residual-constrained formulations.

The 5D seismic synthetic data set is generated using the Compass velocity model provided to us by the BG Group, which has $101\times 101$ receivers spaced by 25m, and $40\times 40$ sources spaced by 25 also with temporal sampling interval of 0.004s. 
It consists of two receiver and two source coordinates, with time as the fifth dimension. We use ocean-bottom nodes style acquisition, where we places the receivers at the ocean-bed and sources are placed at sea-surface. In all experiments, we initialize $L$ and $R$ using standard Gaussian random matrices.  More sophisticated initialization schemes have been proposed~(\cite{jain2013low, slim}), but are expensive for large-scale seismic data since they require computing the largest $r$ singular vectors.

We compare the practical performance (signal-to-noise ratio (SNR) in dB of the reconstructions) and computational time (in seconds) of the new primal-dual formulation to that of LR-BPDN \cite{{slim}}. For simplicity, we perform interpolation for a missing-sources scenario, though the missing-receivers scenario is similar. During seismic data acquisition, the sampling operator $P_{\Omega}$ only acts along source and receiver coordinates and does not depend on the time axis.  We first perform the Fourier transform along the time-axis, and then interpolate each monochromatic data matrix independently. This approach avoids repeated application of Fourier and inverse Fourier transforms during the optimization process.  Once interpolation is finished for each monochromatic data slice, we perform a single inverse Fourier transform along the frequency axis to get the final reconstructed seismic data volume in the time-domain.

We use jittered sampling scheme~\cite{hennenfent2008GEOPsdw} to remove $80\%$ of the sources resulting in the minimum spatial sampling of 25 m and maximum spatial sampling of  225 m between two consecutive sources. The spatial sampling of final interpolated grid is 25 m. Since seismic data is band-limited because of the band-limited nature of the source wavelet, we perform the interpolation within the frequency bandwidth of seismic data, which ranges from 3-70 Hz.  In practice, successful recovery means that an actionable seismic volume (preserving all the coherent energy including late arrivals) can be obtained by computational techniques from merely $20\%$ of acquired data. 

To choose these rank of the factors $L$ and $R$ for interpolation, we consider the 5 Hz and 70 Hz frequency slices,  subsample columns, and perform reconstruction. The best $r$ values (according to SNR value)  for 5 and 70 Hz monocromatic slices were 30 and 100, respectively. Hence, we work with all of the monochromatic slices and adjust the rank linearly from 30 to 100 when moving from lower to higher frequencies. To select 
the target value of the data fitting constraint $\eta$, we use the value $0.03\|b\|_F$, where $b$ 
is the observed data. As described in Section~\ref{sec:relax}, we start Algorithm~\ref{alg:alg1} 
with a relatively loose value of $\eta_0 = \|b\|_F$, and decrease geometrically until we arrive at the target value of $0.03\|b\|_F$.

\begin{figure}[htp]
\centering
\subfloat[]{
\includegraphics[width=.4\textwidth]{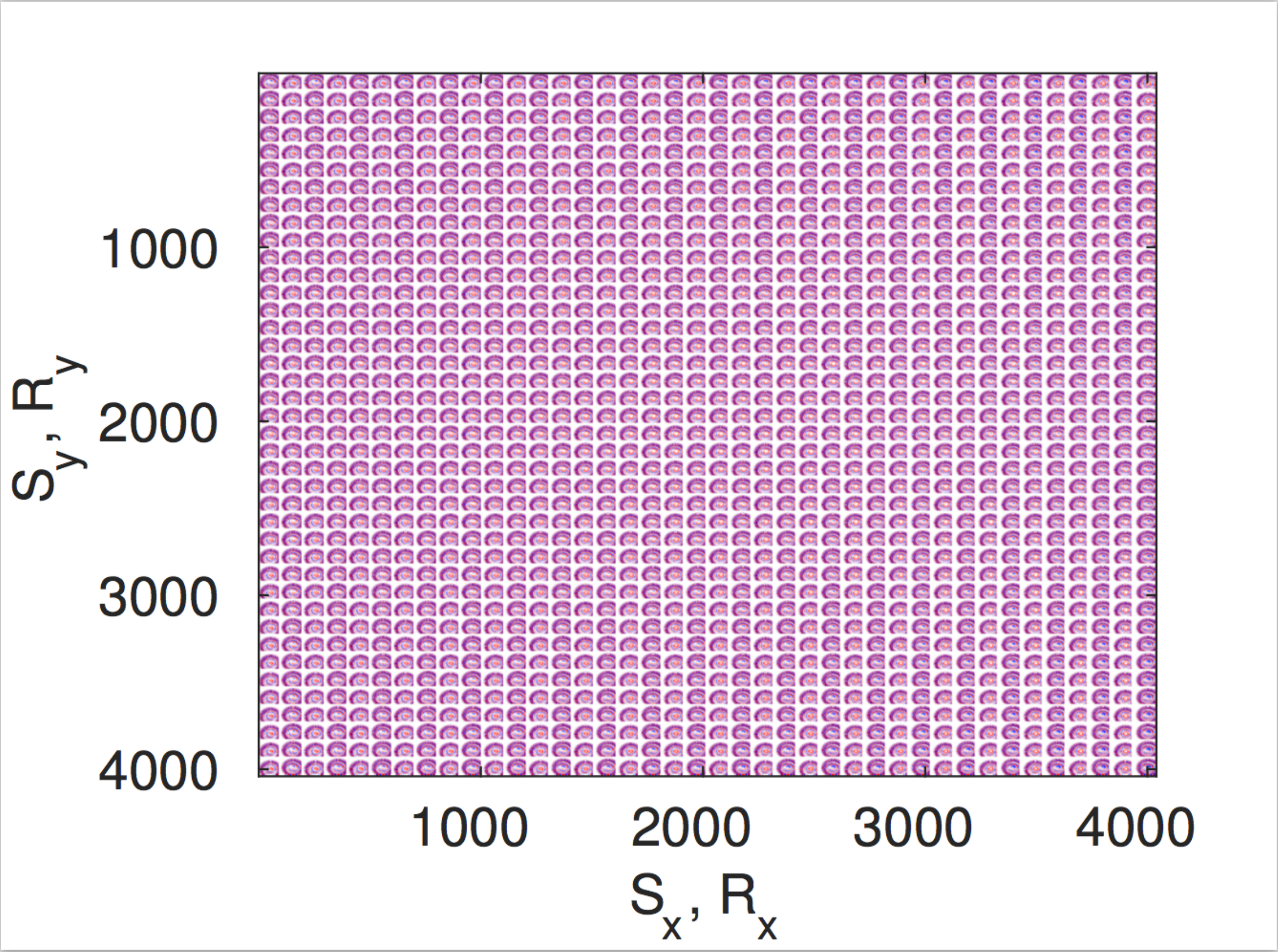}}~
\subfloat[]{
\includegraphics[width=.38\textwidth]{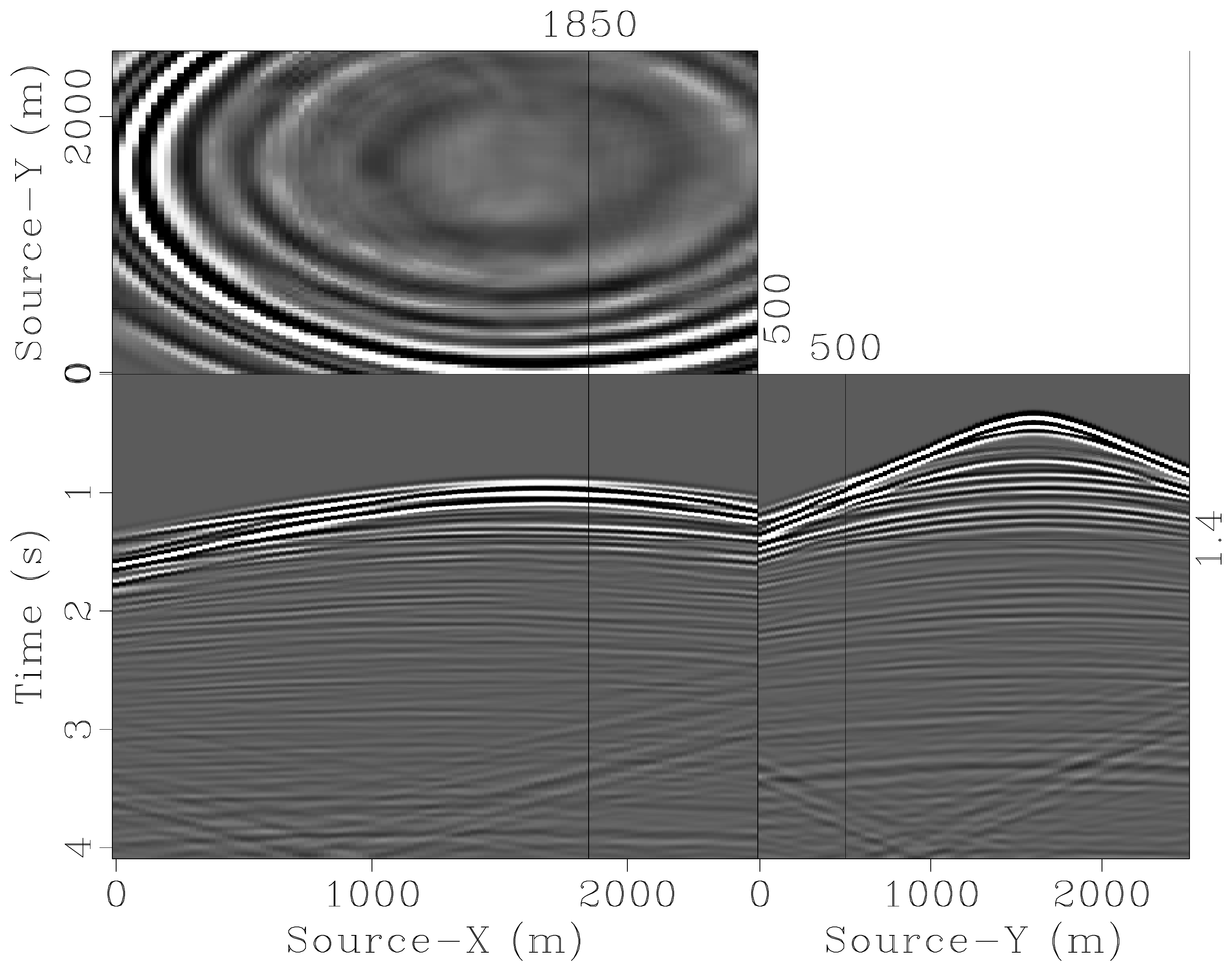}}~\\
\subfloat[]{
\includegraphics[width=.4\textwidth]{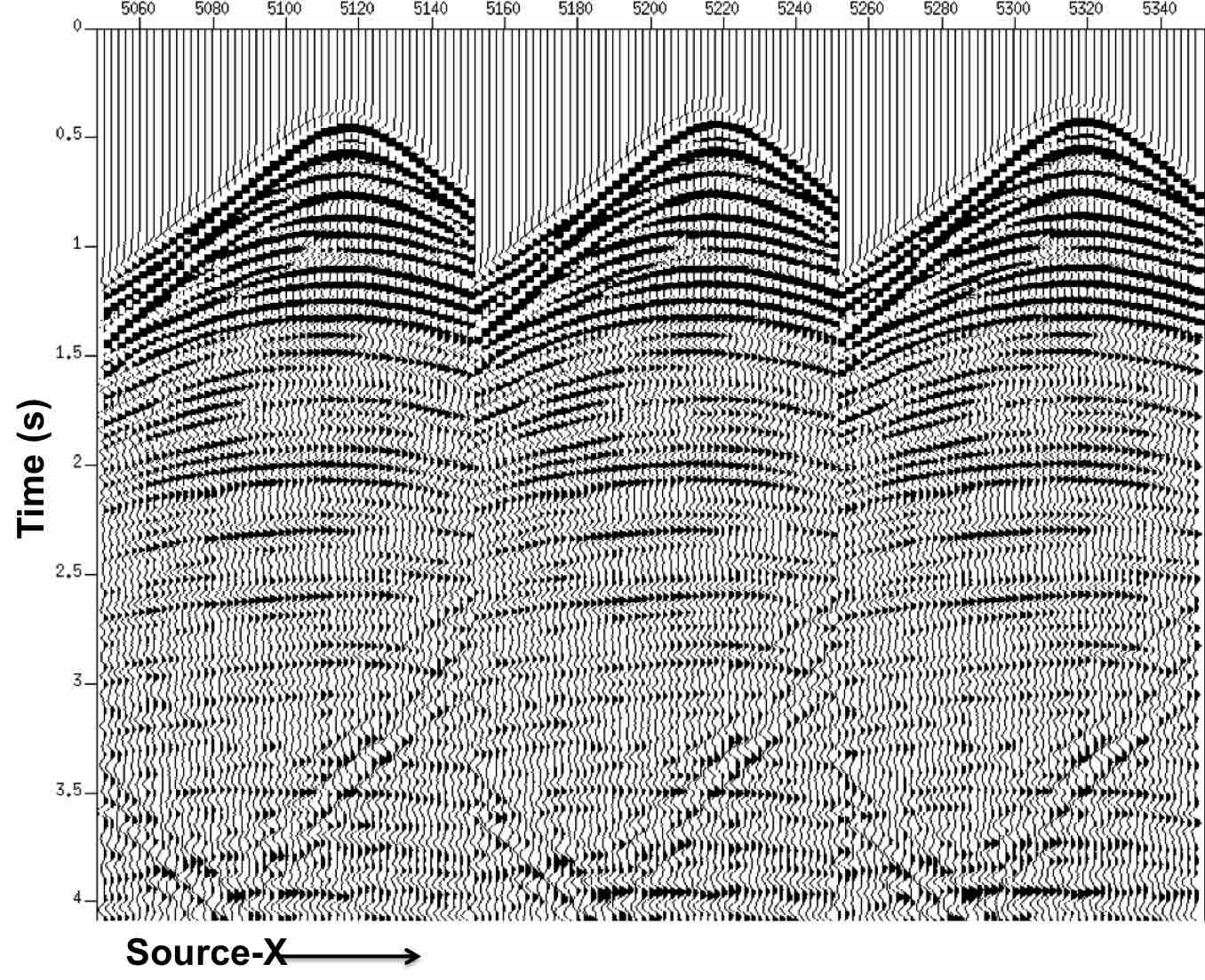}}~
\caption{To show the scale of the fully sampled 5D seismic data volume, we use three levels of granularity. (a) The first step shows a single 4D time-slice, consisting $101 \times 101 \times 40 \times 40$ time samples, from the 5D tensor, where we use ($s_x, r_x$) matricization of the 4D tensor to display it as a matrix. The fully sampled data has 1024 time-slices. (b) In the second step, we extract {\bf one} common-receiver gather from (a) and unfold it along the time-source-x-source-y dimensions. Each common-receiver gather consists of $1024 \times 101 \times 101$ samples. (c) In the final step, for detailed visualization we extract three consecutive columns from (b) and unfold them along the source-axis. 
Interpolation is performed on the entire data volume, comprising multiple 4D monochromatic slices, each of containing $101 \times 101 \times 40 \times 40$ samples.}
\label{fig:3Dtruetime}
\end{figure}
\begin{figure}[htp]
\centering
\subfloat[]{
\includegraphics[width=.4\textwidth]{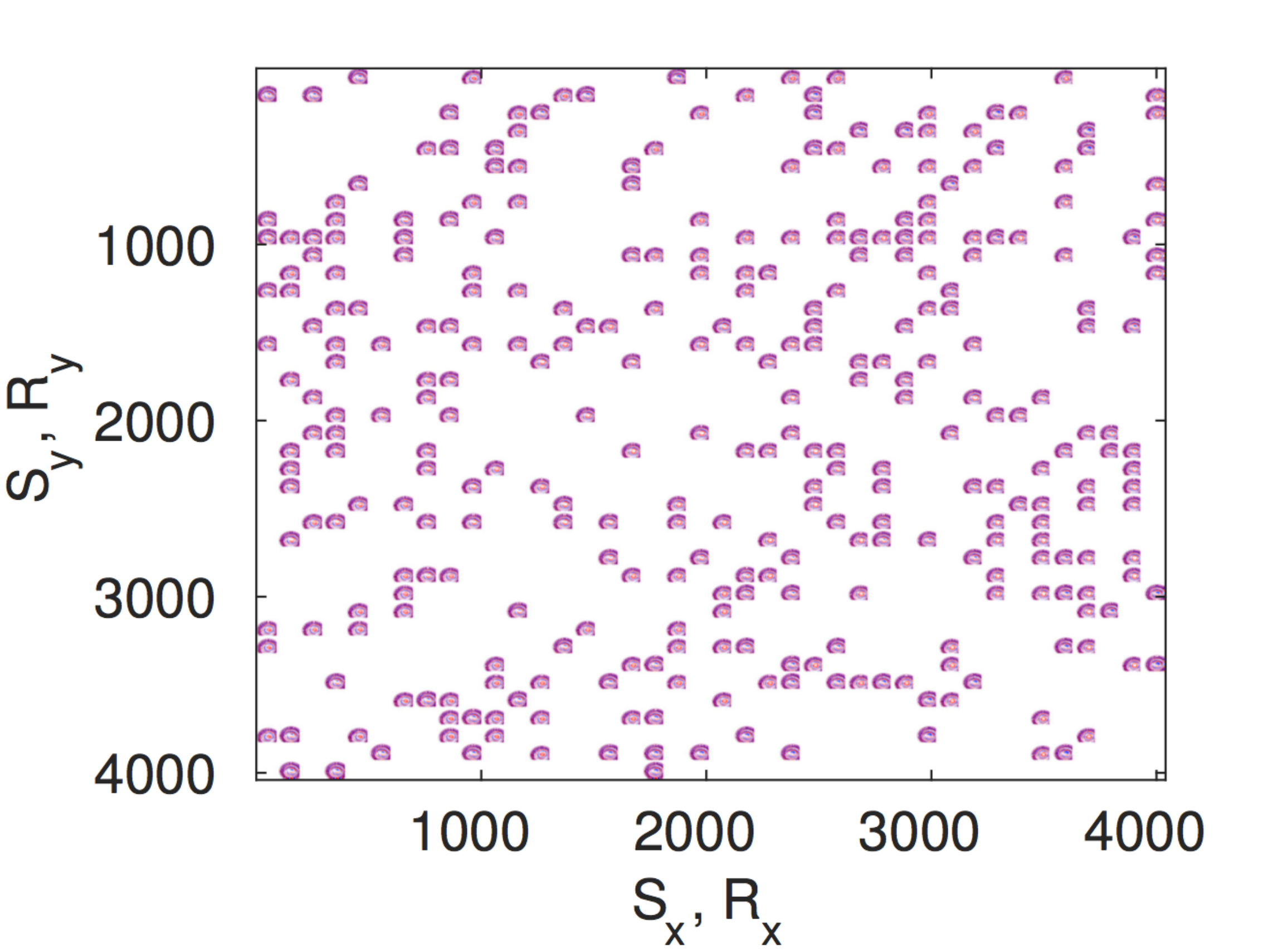}}~
\subfloat[]{
\includegraphics[width=.38\textwidth]{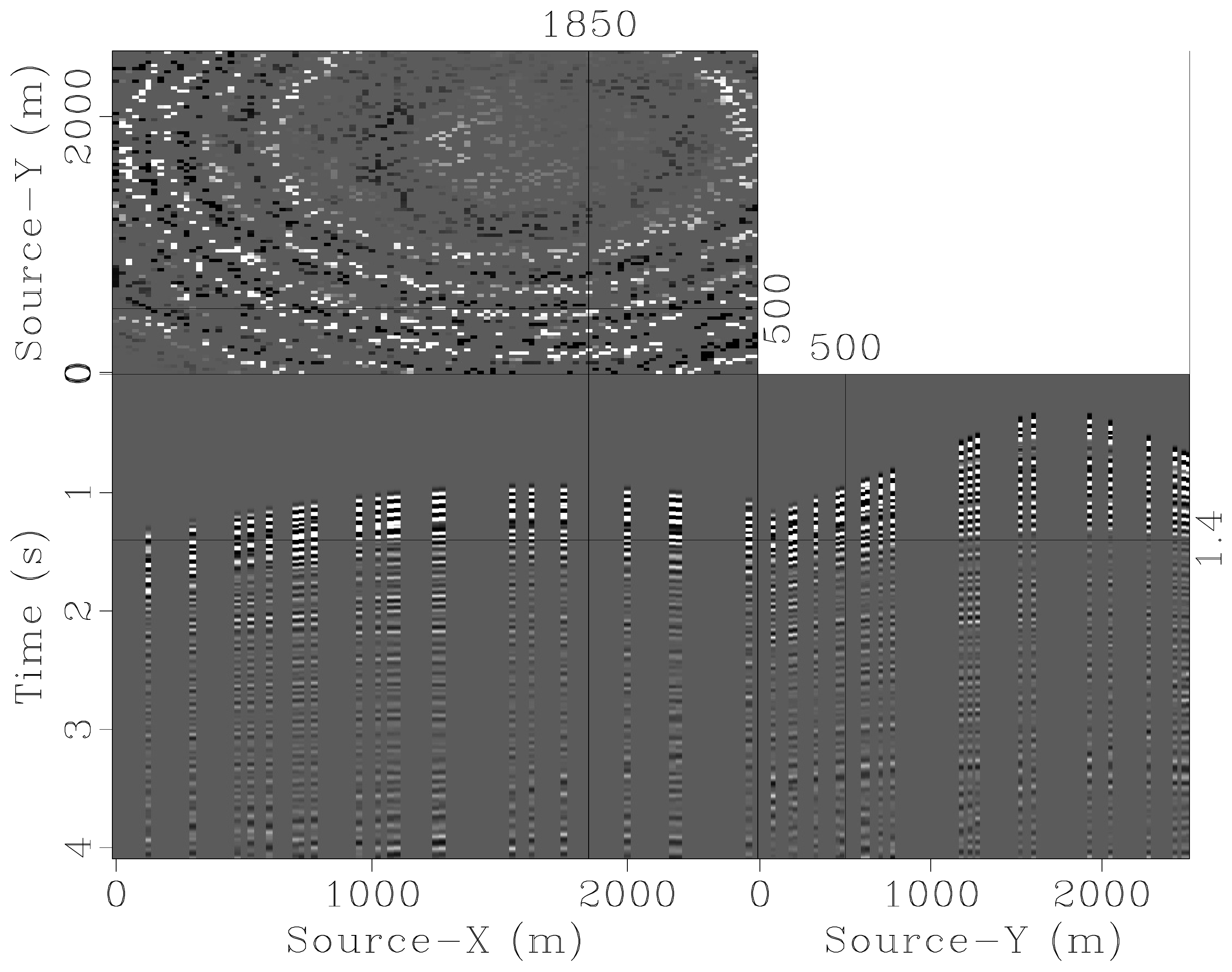}}~\\
\subfloat[]{
\includegraphics[width=.4\textwidth]{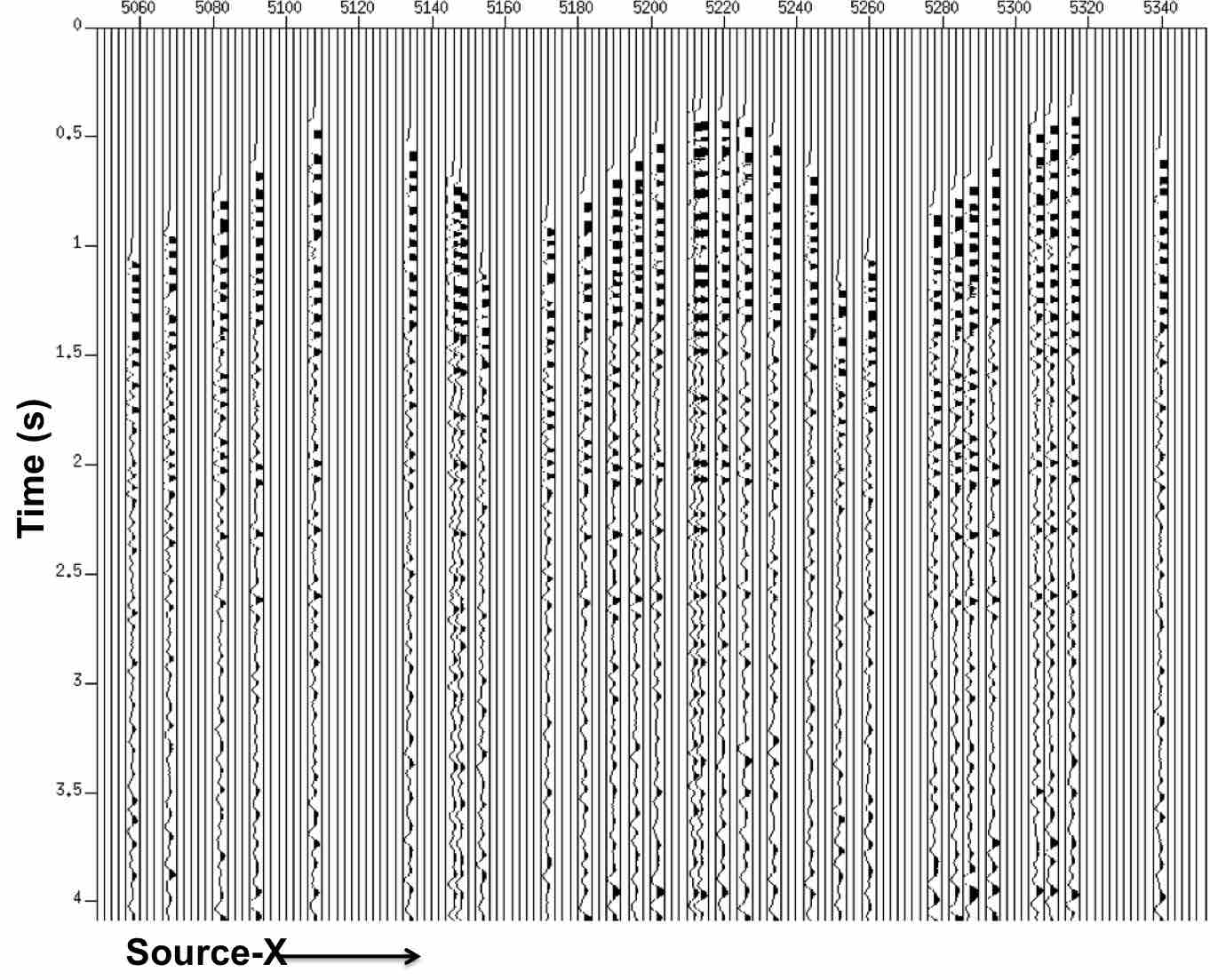}}\\
\caption{Same as Figure~\ref{fig:3Dtruetime} but after removing  80$\%$ jittered sources. Our goal is to recover fully sampled data as shown in Figure~\ref{fig:3Dtruetime} from the 80$\%$ subsampled data.}
\label{fig:3Dtruefreq}
\end{figure}
\begin{figure}[htp]
\centering
\subfloat[]{
\includegraphics[width=.4\textwidth]{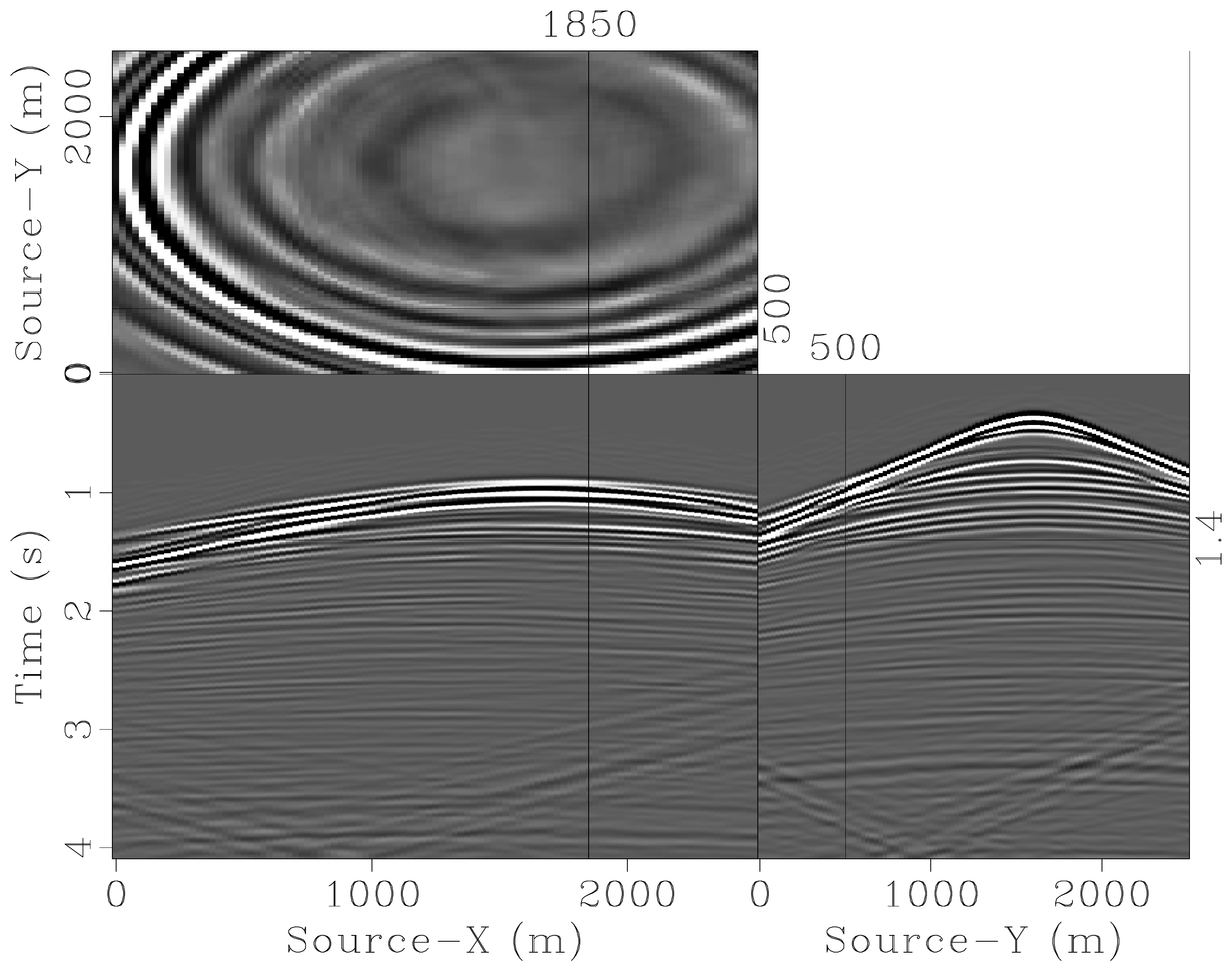}}~
\subfloat[]{
\includegraphics[width=.4\textwidth]{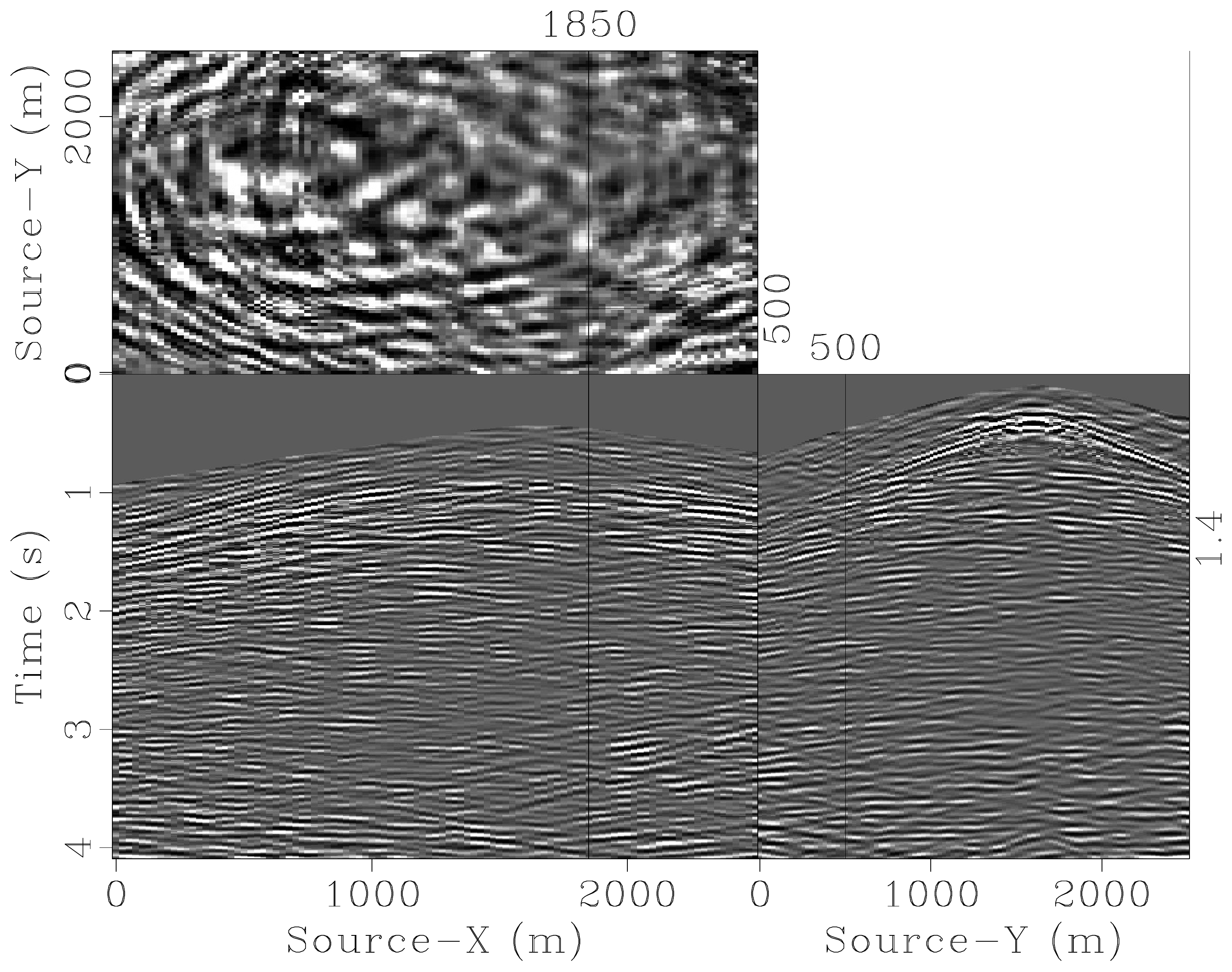}}\\
\subfloat[]{
\includegraphics[width=.4\textwidth]{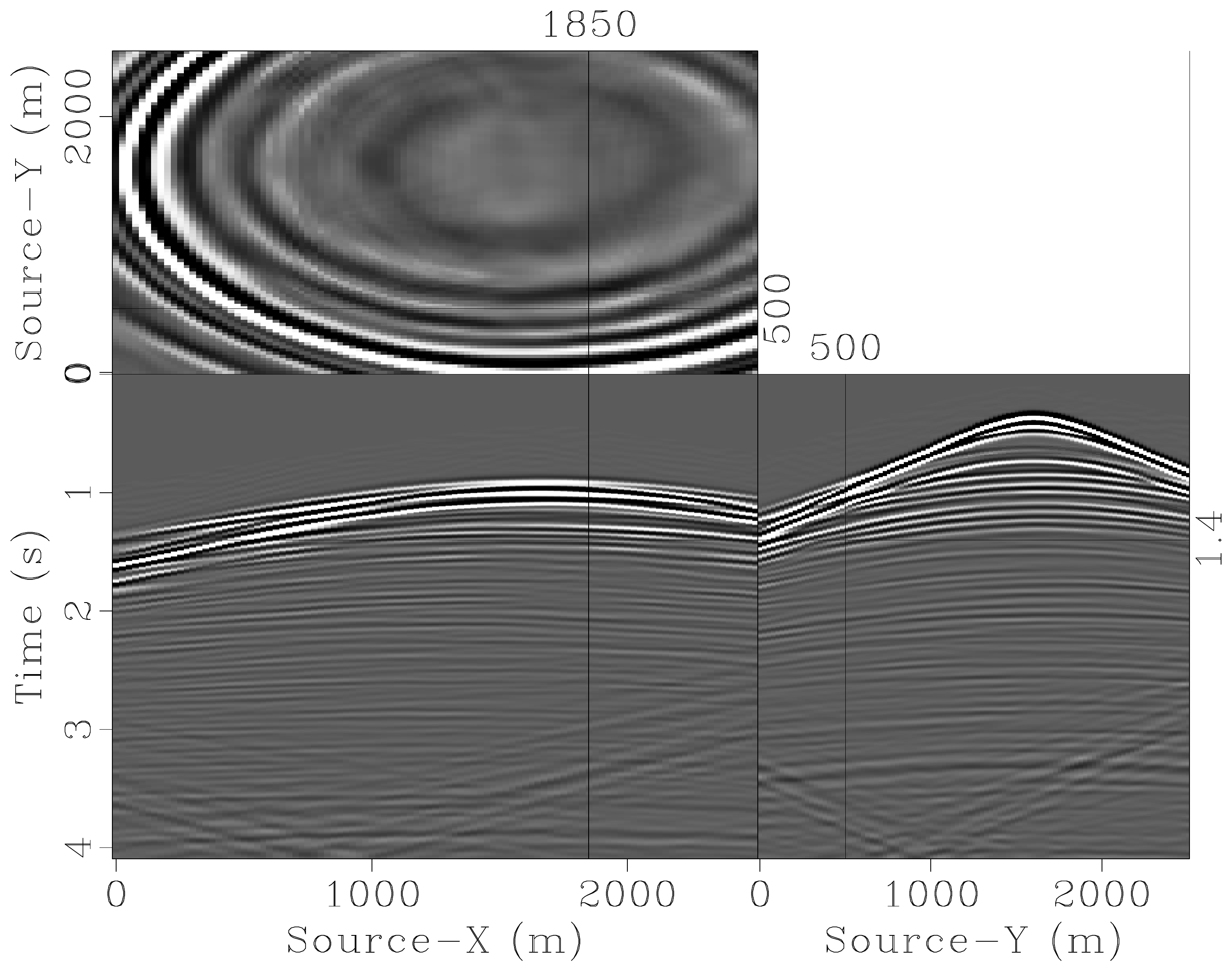}}~
\subfloat[]{
\includegraphics[width=.4\textwidth]{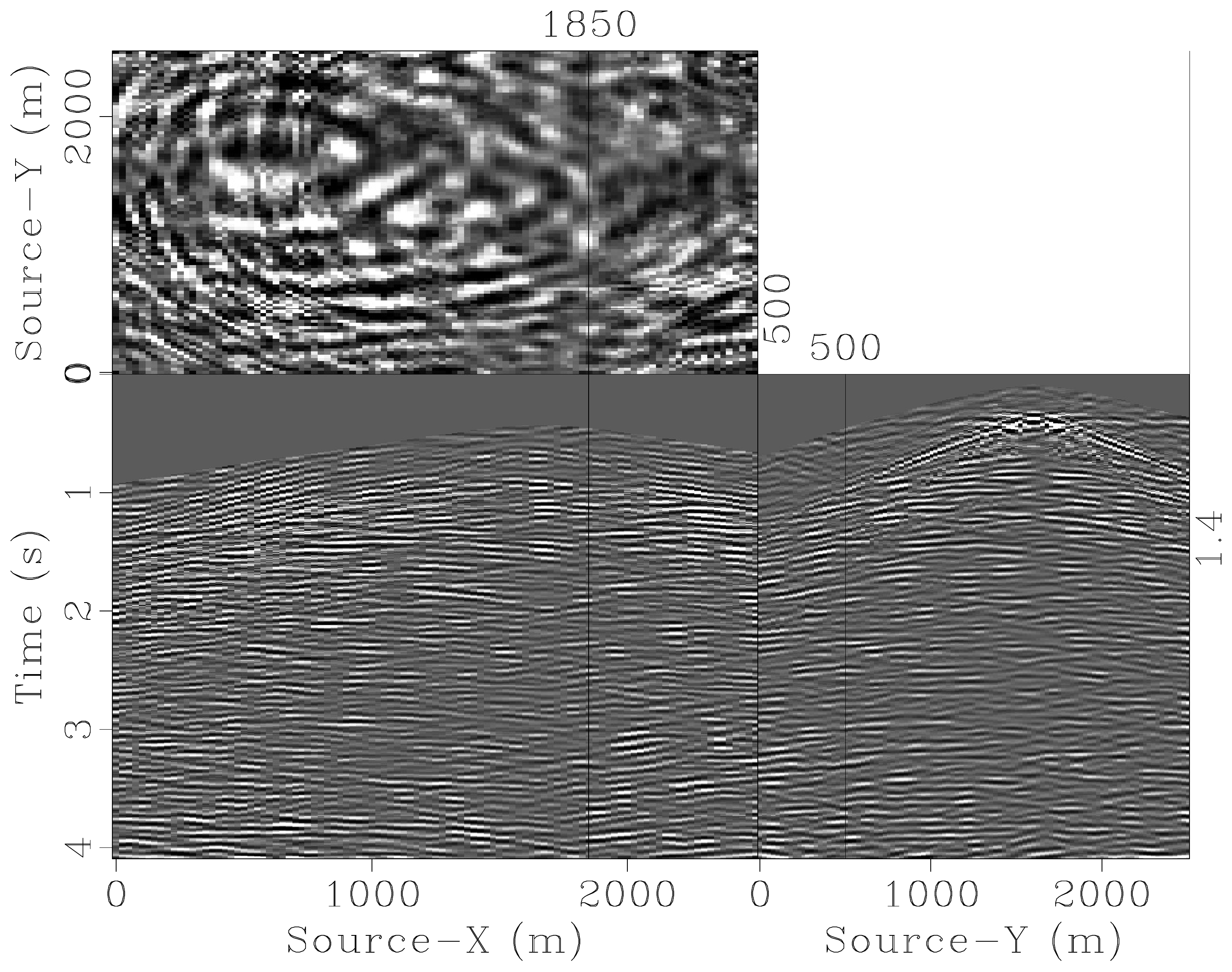}}
\caption{A common-receiver gather from interpolated data and 100-times amplified residuals using (a,b) primal-dual splitting and (c,d) using level-set method for 80$\%$ missing sources scenario. We can clearly see that we recover most of the coherent energy using both the methods. We amplify the amplitude of residual by a factor of 100 to show that we are not loosing any coherent energy.}
\label{fig:3D80}
\end{figure}
\begin{figure}[htp]
\centering
\subfloat[]{
\includegraphics[width=.4\textwidth]{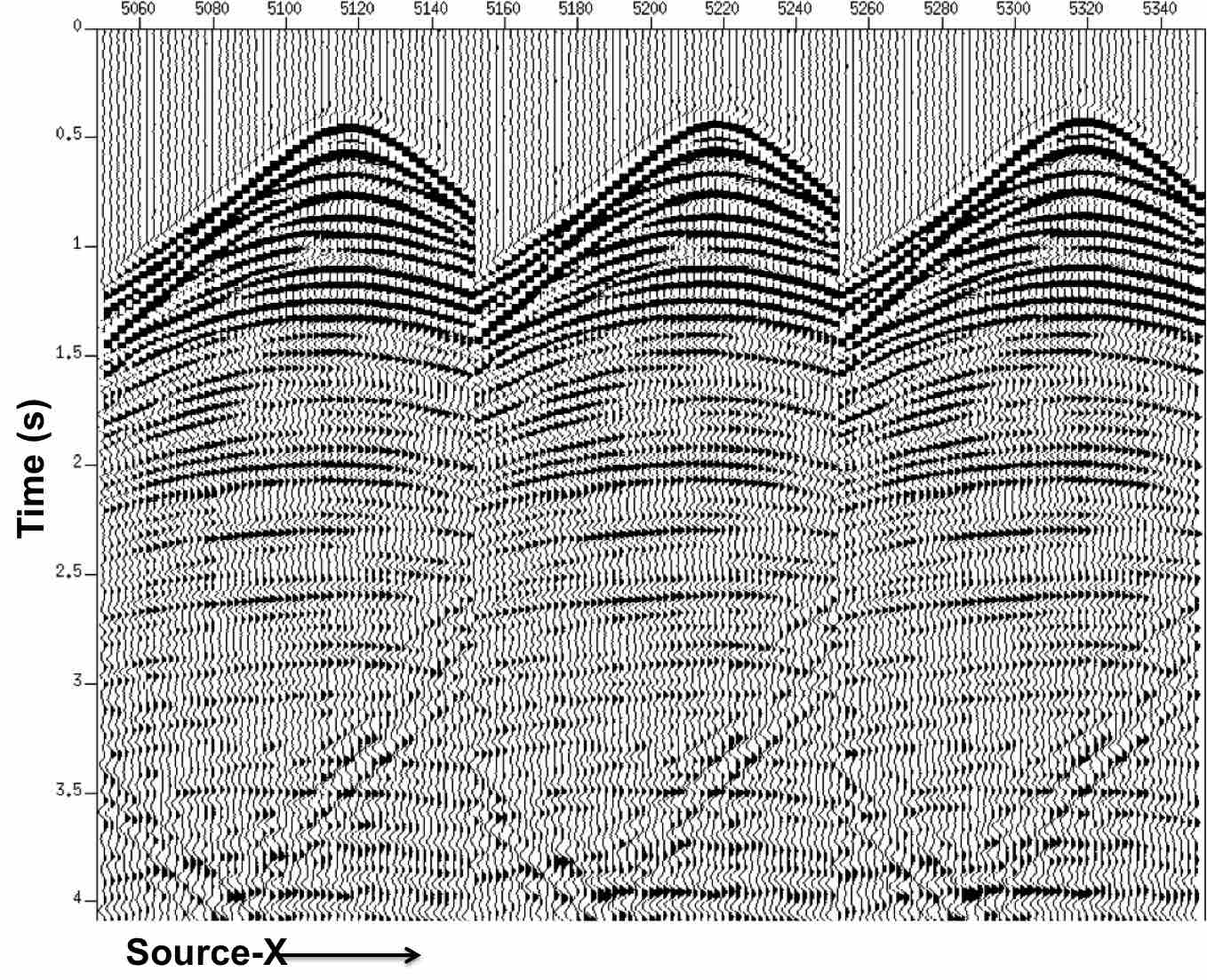}}~
\subfloat[]{
\includegraphics[width=.4\textwidth]{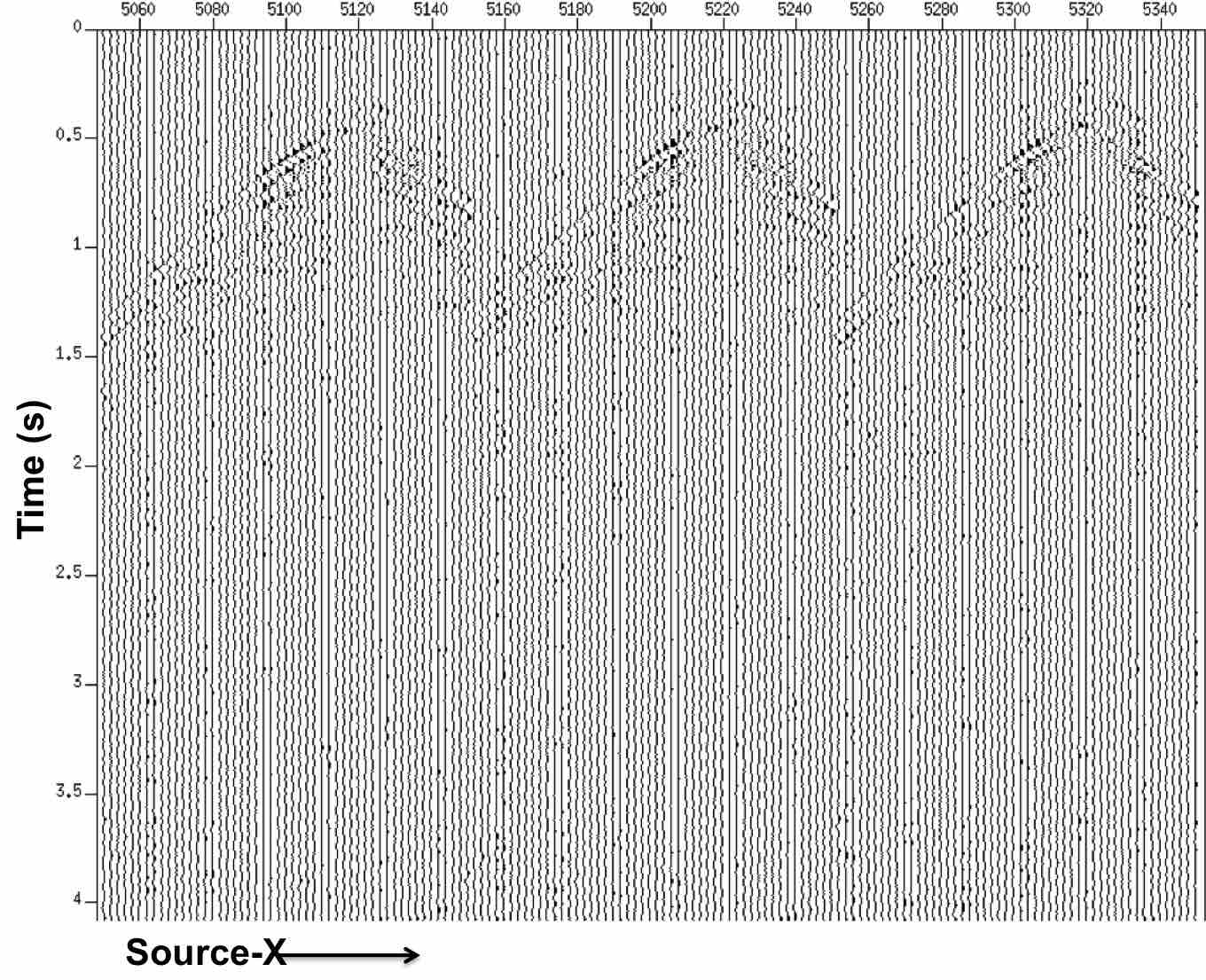}}\\
\subfloat[]{
\includegraphics[width=.4\textwidth]{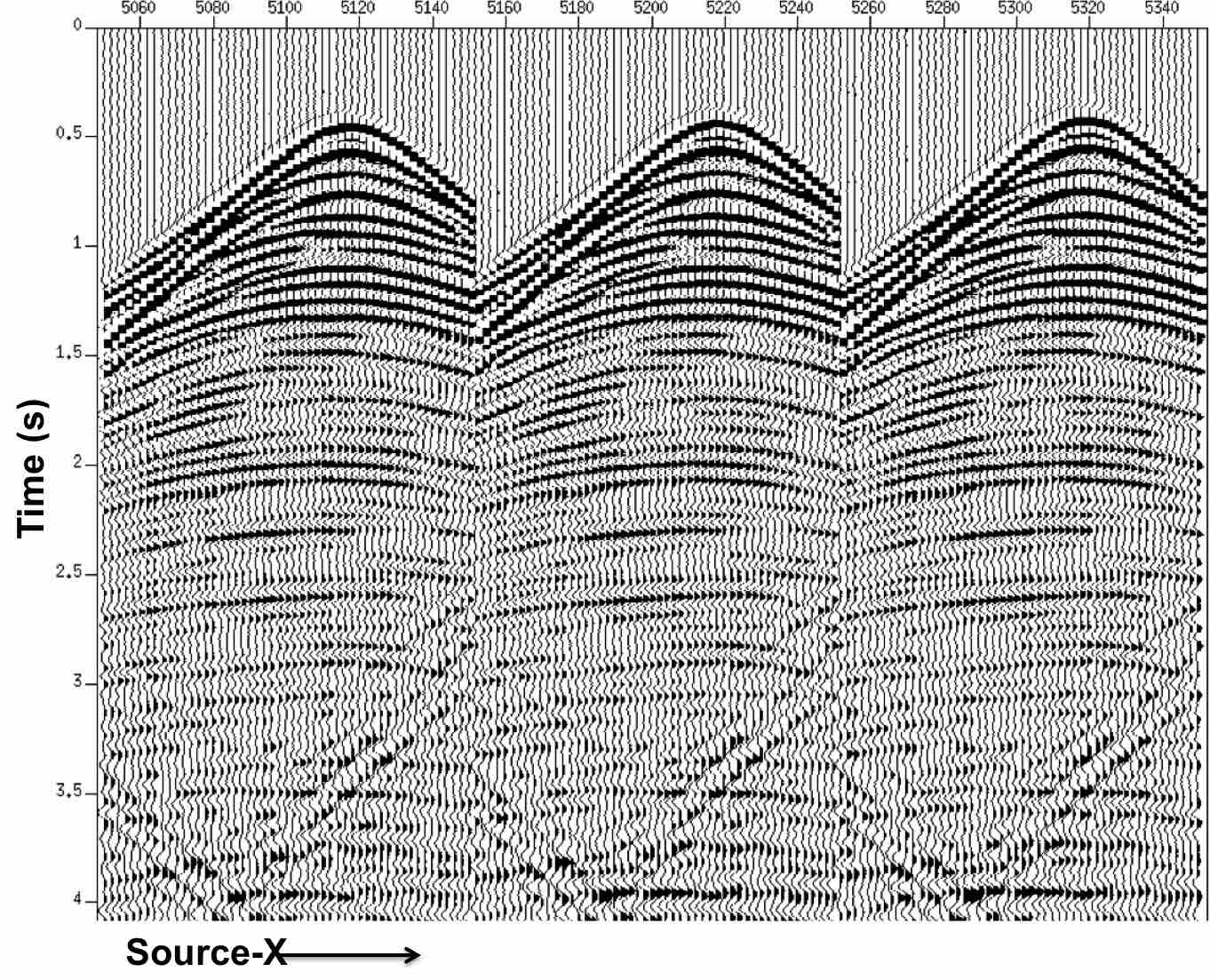}}~
\subfloat[]{
\includegraphics[width=.4\textwidth]{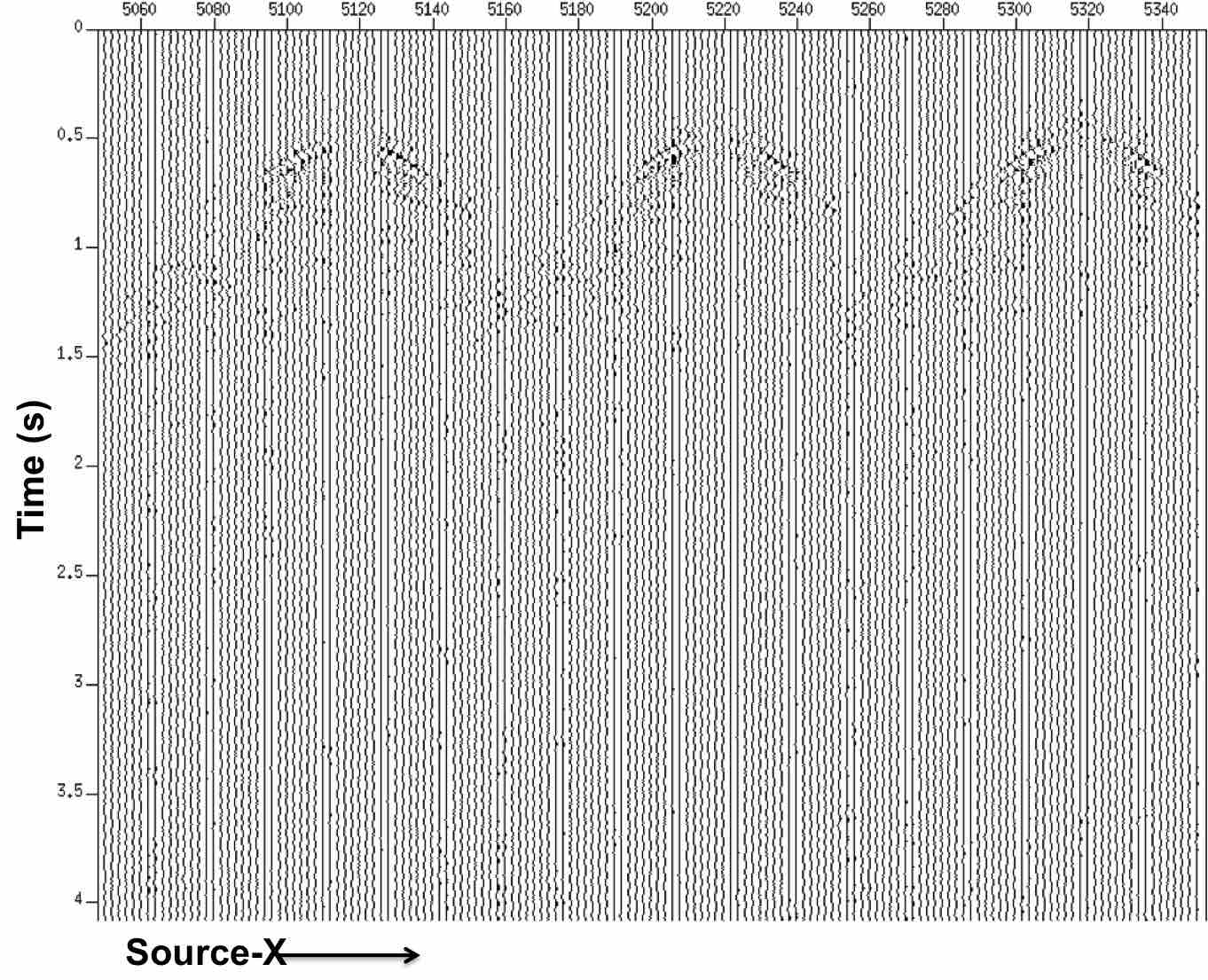}}
\caption{Subsection of common-receiver gather (extracted from Figure~\ref{fig:3D80}) after interpolation and corresponding residual using, (a,b) the proposed method, and (c,d) level-set method. Here, we keep the amplitude scale of reconstruction and residual same to show the efficacy of proposed formulation.}
\label{fig:3D80timewiggle}
\end{figure}

Figures \ref{fig:3Dtruetime} and \ref{fig:3Dtruefreq} show fully sampled and subsampled (80\% missing sources) data from the 5D seismic volume. Figures \ref{fig:3Dtruetime} (c) and \ref{fig:3Dtruefreq} (c) show a zoom section of common-receiver gather from the true and subsampled data. The reconstructed data volume and corresponding residual plots are shown in Figures \ref{fig:3D80} (a, b) and (c,d) using the proposed alternating primal-dual splitting method and the LR-BPDN based level-set approach. We boost the amplitudes of residual plots (Figures \ref{fig:3D80} (c,d)) by a factor of 100 to show that we are not loosing any coherent energy. Both of the methods are able to reconstruct the coherent energy of 5D seismic data volumes while preserving the late-arrival energy, which can be also seen from the residual plots. We also plot the zoom sections (Figure \ref{fig:3D80timewiggle}) of reconstructed data and corresponding residual to strengthen our observations. 

\begin{figure}[htp]
\centering
\subfloat[]{
\includegraphics[width=.6\textwidth]{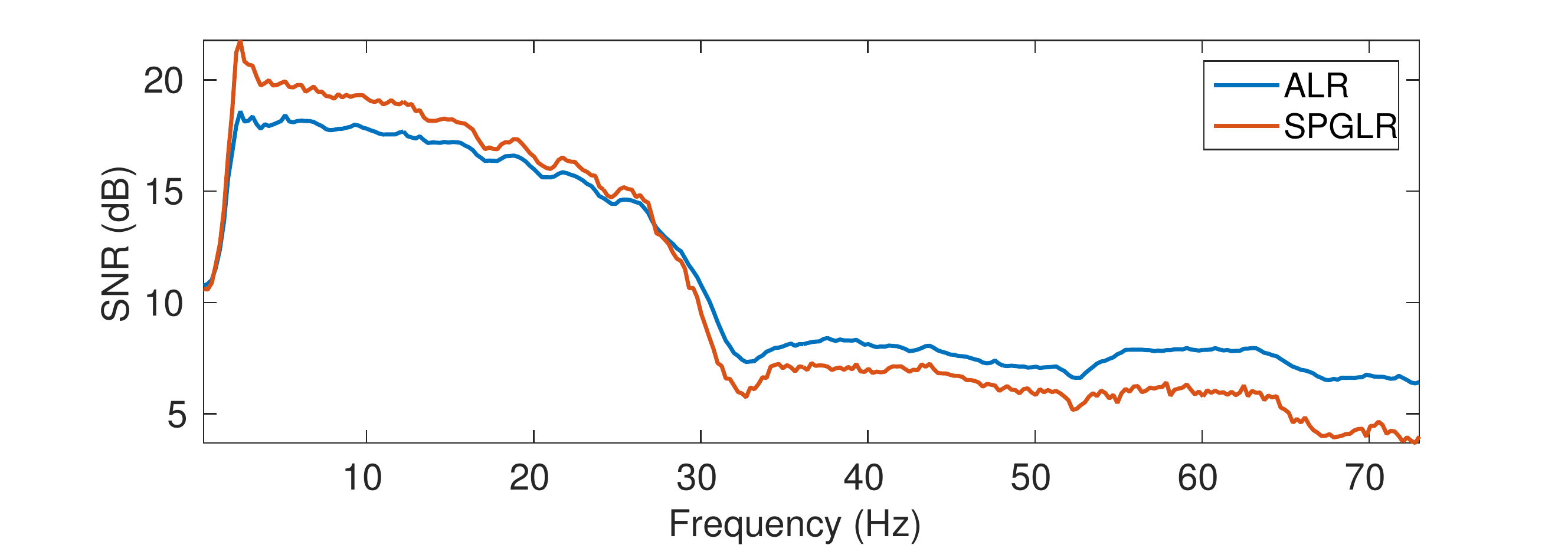}}\\
\subfloat[]{
\includegraphics[width=.6\textwidth]{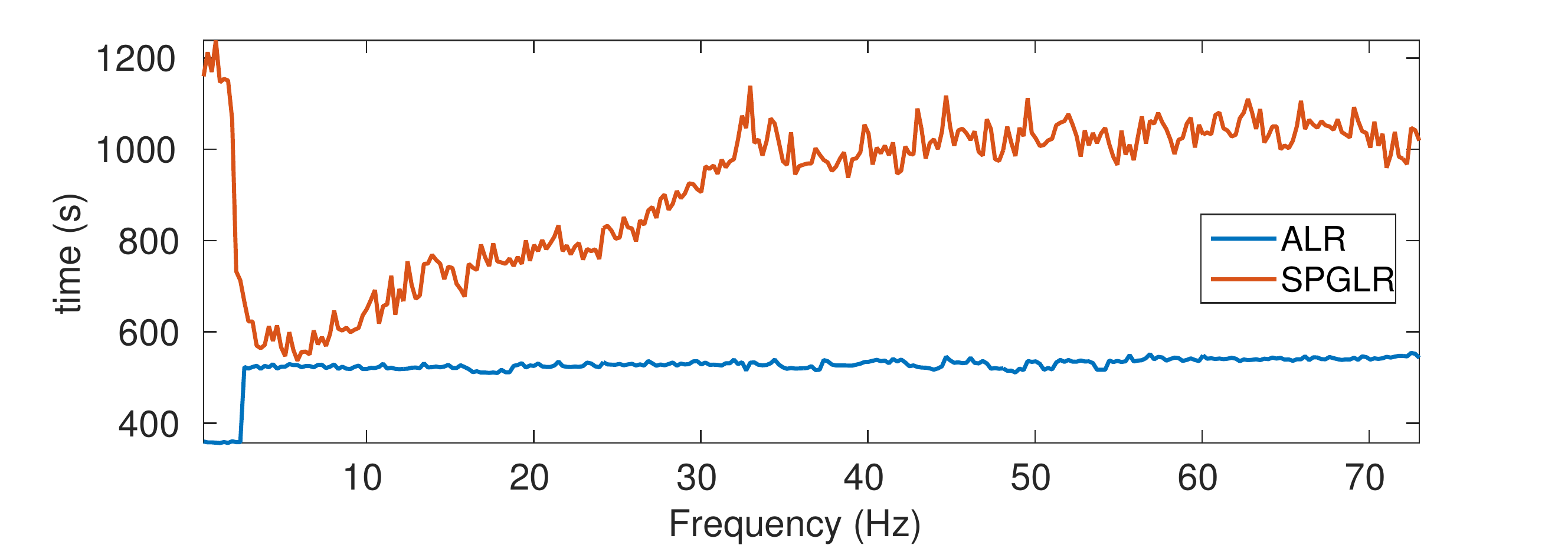}}
\caption{Qualitative comparison of level-set method and the proposed primal-dual methodology over the frequency spectrum of 5D synthetic seismic data volume generated using BG Compass model for $80\%$ missing sources. (a) signal-to-noise ratio (dB), (b) computational time (s).}
\label{fig:SNR803D}
\end{figure}

We compare the signal-to-noise ratio (SNR in dB) and computational times (in sec.) across the frequency spectrum as shown in Figures \ref{fig:SNR803D} (a, b). We can see that the alternating primal-dual method works better along the higher frequency regime whereas the level-set method works better along the lower-frequency regime. Computationally, the primal-dual method is faster by a factor of two. In addition, the computational 
time for the primal-dual scheme is remarkably consistent across the frequency spectrum, while LR-BPDN requires more time as frequency increases. We run both the methods on a system with two 10-core Intel E5-2690 v2 Ivy Bridge CPUs at 3.00GHz.

%
%


\section{Discussion and conclusions}
Residual-constrained formulations are well suited for practical interpolation schemes, 
since they provide a simple way for practitioners to fit to a target data-fitting level. 
In this work, we propose a new approach for residual-constrained formulations, 
using a block coordinate alternating optimization scheme. 
As with other factorized schemes, in the current approach 
the full data volume is represented implicitly via the outer product $X=LR^T$, 
and achieves time and memory savings by working with the factors instead of the full volume $X$. 

Even though the overall problem is nonconvex, each problem in $L$ and $R$ is a convex 
residual-constrained problem, and we solve it using a primal-dual splitting approach detailed in Algorithm~\ref{alg:alg2}. The approach requires only matrix-vector products with $L, R$ and their adjoints. 
It is simpler than the level-set approach, which requires 
variational computations of the value function. Experimental results on large-scale 
seismic interpolation show that the proposed approach 
can match the practical performance of level set methods, 
obtaining comparable results in half the time, and working consistently across the 
range of frequencies of interest. 

An important consequence of this work is that it opens future directions in regularized data interpolation. 
One can for example consider constraints on the factors $L, R$ themselves, or additional constraints 
on $LR^{H}$. When working with additional convex constraints, Algorithm~\ref{alg:alg1} can be left completely unchanged, and  Algorithm~\ref{alg:alg2} can be adapted. The key point is that with $L$ or $R$ fixed, the subproblems remain convex. We leave these developments for future work.

\label{sec:discussion}

\section*{Acknowledgements}
We would like to thank the BG Group for permission to use the synthetic Compass velocity model. The authors wish to acknowledge the SENAI CIMATEC Supercomputing Center for Industrial Innovation, with support from BG Brazil and the Brazilian Authority for Oil, Gas and Biofuels (ANP), for the provision and operation of computational facilities and the commitment to invest in Research \& Development. This work was financially supported in part by the Natural Sciences and Engineering Research Council of Canada Collaborative Research and Development Grant DNOISE II (CDRP J 375142-08). This research was carried out as part of the SINBAD II project with the support of the member organizations of the SINBAD Consortium. Aleksandr Aravkin was partially supported by the Washington Research Foundation Data Science Professorship.  

\bibliography{siam}
\bibliographystyle{abbrv}

\end{document}